\documentclass[english]{amsart}
\usepackage[T1]{fontenc}
\usepackage[utf8]{inputenc}
\usepackage{csquotes, xpatch}
\usepackage{microtype}
\usepackage{amsmath}
\usepackage{amssymb}
\usepackage{amsthm}
\usepackage{mathrsfs}
\usepackage{mathtools}
\usepackage[all]{xy}
\usepackage{caption}
\usepackage{subfig}
\usepackage{booktabs}
\usepackage{graphicx}
\usepackage{tikz}
\usetikzlibrary{automata,positioning,calc,intersections,through,backgrounds,patterns,fit,external}
\usepackage[style=alphabetic, maxnames=5, backend=biber,sorting=nyt]{biblatex}
	\renewbibmacro{in:}{\ifentrytype{article}{}{\printtext{\bibstring{in}\intitlepunct}}}
\usepackage{url}
\usepackage[colorlinks]{hyperref}
	\hypersetup{hidelinks}		
\usepackage{bookmark}
\usepackage[nameinlink,capitalise,noabbrev]{cleveref}

\theoremstyle{plain}
\newtheorem{theorem}{Theorem}
\newtheorem{lemma}[theorem]{Lemma}
\newtheorem{proposition}[theorem]{Proposition}
\newtheorem{corollary}[theorem]{Corollary}
\newtheorem{conjecture}[theorem]{Conjecture}
\theoremstyle{definition}
\newtheorem{definition}[theorem]{Definition}
\theoremstyle{remark}
\newtheorem*{remark}{Remark}

\newcommand{\numberset}{\mathbb}

\newcommand{\NN}{\numberset{N}}
\newcommand{\PP}{\numberset{P}}

\newcommand{\RR}{\numberset{R}}

\newcommand{\cA}{\mathcal{A}}

\newcommand{\cD}{\mathcal{D}}
\newcommand{\cE}{\mathcal{E}}
\newcommand{\cH}{\mathcal{H}}
\newcommand{\cI}{\mathcal{I}}
\newcommand{\cL}{\mathcal{L}}

\renewcommand{\epsilon}{\varepsilon}

\newcommand{\alphabet}{\mathbf{A}} 
\newcommand{\rauzy}{\mathcal{R}}
\newcommand{\rauzygraph}{G}

\newcommand{\gasket}{\mathscr{G}}
\renewcommand{\sl}{\mathfrak{sl}}
\newcommand{\bX}{\mathbf{X}}
\newcommand{\sW}{\mathscr{W}}

\newcommand*{\transpose}[2][-7mu]{\ensuremath{\mskip1mu\prescript{\smash{\top\mkern#1}}{}{\mathstrut#2}}}%

\DeclareMathOperator{\GL}{GL}
\DeclareMathOperator{\SL}{SL}
\DeclareMathOperator{\SO}{SO}

\DeclareMathOperator{\diam}{diam}
\DeclareMathOperator{\area}{area}

\begin{document}

\title[Renormalization for Bruin-Troubetzkoy ITMs]{Renormalization for Bruin-Troubetzkoy ITMs}
\date{\today}

\author{Mauro Artigiani}
\address{Universidad Nacional de Colombia - Bogotá}
\email[M.~Artigiani]{martigiani@unal.edu.co}

\author{Pascal Hubert}
\address{Aix-Marseille Universit\'e, CNRS\\ Institut de Math\'ematiques de Marseille, I2L, 13453, Marseille, France}
\email[P.~Hubert]{pascal.hubert@univ-amu.fr}

\author{Alexandra Skripchenko}
\address{Faculty of Mathematics\\
	National Research University Higher School of Economics\\
	Usacheva St. 6\\
	119048 Moscow, Russia}
\email[A.~Skripchenko]{sashaskrip@gmail.com}



\begin{abstract}
	We study a class of interval translation mappings introduced by Bruin and
	Troubetzkoy, describing a new renormalization scheme, inspired by the
	classical Rauzy induction, for this class. We construct a measure, invariant
	under the renormalization, supported on the parameters yielding infinite
	type interval translation mappings in this class. With respect to this
	measure, a.e.\ transformation is uniquely ergodic. We show that this set has
	Hausdorff dimension between $1.5$ and $2$, and that the Hausdorff dimension
	coincides with the affinity dimension. Finally, seeing our renormalization
	as a multidimensional continued fraction algorithm, we show that it has
	almost always the Pisot property.
	
	We discover an interesting phenomenon: the dynamics of this class of
	transformations is often (conjecturally: almost always) weak mixing, while
	the renormalizing algorithm typically has the Pisot property.
\end{abstract}

\maketitle

\section{Introduction}
This paper is focused on the ergodic properties of two classes of related
dynamical systems. The first class we are interested in is a particular family
of interval translation mappings (or ITMs, for short), the second one is a
Markovian multidimensional continued fraction algorithm (MCF).

ITMs were introduced in~\cite{BK} as a natural generalization of \emph{interval
exchange transformations} (IETs). IETs and their ergodic properties were widely
studied in the last decades, see, e.g.,~\cites{Viana:IETs, Yoccoz:IEMs} and the
references therein. Typical IETs are known to be uniquely
ergodic~\cites{Masur:ue, Veech:ue} and weakly mixing~\cite{AvilaForni}, while
ergodic properties of certain special classes of IETs can be remarkably
different (for example, Arnoux-Rauzy IETs are almost never minimal and those who
are minimal are typically not weakly mixing~\cite{ACFH}). All these results were
achieved by the study of the properties of the renormalization algorithm called
\emph{Rauzy induction} (and variations of it). This algorithm can be seen as a
representative of the class of Markovian multidimensional continued fraction
algorithms. 

\begin{figure}[tb]
	\centering
	\includegraphics[width=.6\textwidth]{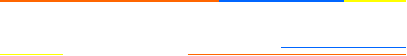}
	\caption{An example of a Bruin-Troubetzkoy ITM. The intervals below are images of the ones above, color coded.}
	\label{fig:BT_ITM}
\end{figure}

The key difference between IETs and ITMs is that the latter are not necessarily
surjective: the images of the intervals do not need to form a partition, they
simply form a collection of subintervals of the original interval,
see~\cref{fig:BT_ITM} for an example. More formally,

\begin{definition}
	An \emph{interval translation map} is a piecewise translation map $T$
	defined on a half-open interval $I = [0,a) \subset \RR$, for some $a>0$,
	with values in $I$. We call $T$ a $n$-interval translation map (or $n$-ITM)
	if $I$ has $n$ maximal half-open sub-intervals to which the restriction of
	$T$ is a translation. The endpoints of these intervals are called
	\emph{singularities} of the map.
\end{definition}

It was noticed already in~\cite{BK} that each ITM is either of \emph{finite} or
\emph{infinite} type. This classification is based on the properties of the
\emph{attractor} of the ITM. Namely, for a given mapping $T$ we consider the
sequence $\Omega_n = I\cap TI\cap T^2I\cdots\cap T^{n}I$. If this sequence
stabilizes for some $N\in\NN$, i.e., $\Omega_k = \Omega_{k+1}$ for all $k\ge N$,
then the ITM $T$ is of \emph{finite type}. If there is no such $N$ and
$\overline{\Omega}$, the closure of the limit set $\Omega =I\cap TI\cap
T^2I\cdots$, is a Cantor set, then the ITM is of \emph{infinite type}, see
also~\cite{ST}.

Dynamics of ITMs of finite type basically coincides with the one of IETs.
However, ITMs of infinite type are remarkably different. M.~Boshernitzan and
I.~Kornfeld described the first example of ITM of infinite type. In the same
paper, they formulated the following

\begin{conjecture}\label{conjecture} 
	The set of parameters that give rise to ITMs of infinite type has zero
	Lebesgue measure.
\end{conjecture}

To the best of our knowledge, this conjecture is currently completely open. The
only known cases are for ITMs on $2$ and $3$ intervals, see~\cites{BK, Volk}
respectively, and for a (very special) family of $n$-ITMs, which generalizes the
one we study in this paper to an arbitrarily high number of intervals of
continuity, see~\cite{Br}. Very recently, a \emph{topological} version of the
conjecture has been proven in~\cite{DSvS}.

In this paper, we focus on a special subclass of ITMs, which was defined by
H.~Bruin and S.~Troubetzkoy in~\cite{BT} and which was, historically, the first
concrete example of a family of ITMs. The class is described as follows: let
$U:= \{(\alpha, \beta): 0\le\beta\le\alpha\le1 \} $ and $L:= \{(\alpha, \beta):
0\le\alpha\le\beta+1\le1 \} $ and $R:=U\cup L;$ for the internal point $(\alpha,
\beta)\in U$ we define 
\begin{equation}
T_{\alpha,\beta}(x) = \begin{cases}
                                  x+\alpha, & x\in[0,1-\alpha)\\
                                  x+\beta, & x\in[1-\alpha, 1-\beta)\\
                                  x+\beta-1, & x\in[ 1-\beta, 1),
                                  \end{cases}
\end{equation}
see~\cref{fig:BT_ITM} for an example.

The transformation $T(x) = T_{\alpha,\beta}(x)\colon [0,1)\to[0,1)$ is a
$3$-ITM. By identifying the points $0$ and $1$ we get an interval translation
map on a circle with two intervals of continuity. We remark that the original
example of ITM in the paper by Boshernitzan and Kornfeld also belongs to this
family. In their paper, Bruin and Troubetzkoy proved \cref{conjecture} for this
special family of ITMs (see~\cite[Theorem~6]{BT}). They also showed that,
considering the set of $3$-ITMs $T_{\alpha, \beta}$, the set
$\gasket_{\text{ue}}$ of parameters that gives rise to uniquely ergodic ITMs of
infinite type is a dense $G_\delta$ subset of the set $\gasket$ of parameters
that give rise to ITMs of infinite type (see~\cite[Corollary~13]{BT}). In this
paper, we improve their result in the following way: 

\begin{theorem}\label{thm:uniquelyergodic} There exists a natural measure $\mu$,
	whose support set coincides with $\gasket$, such that for $\mu$-almost all
	$(\alpha, \beta)\in \gasket$, the transformation $T_{\alpha, \beta}$ is
	uniquely ergodic. 
\end{theorem} 

We remark that the result by Bruin and Troubetzkoy was later generalized by
Bruin for a slightly more extended subclass of ITMs (see~\cite[Theorem~1]{Br}).
Also, one can see Bruin-Troubetzkoy ITMs as a special type of double rotations,
which were introduced in~\cite{SIA} and studied in~\cites{BC, AFHS}. 

Bruin and Troubetzkoy obtained their results by describing a special type of
renormalization procedure (a Gauss-like map)  for their class of ITMs. Using
this procedure, they found a symbolic (more precisely, substitutional)
presentation of the interval translation mappings they were interested in, and
used it to prove their Theorem~6 and to construct the non-uniquely ergodic
examples. 

Our strategy is quite different. In fact, we treat Bruin-Troubetzkoy family as a
particular class of \emph{systems of isometries}. First, we introduce a new
renormalization procedure that is based on the induction that I.~Dynnikov
defined for systems of isometries (see \cite{D}). Our renormalization algorithm
is a projectivization of the linear map defined by the induction procedure and
can be seen as a Markovian multidimensional continued fraction (MCF) algorithm.
Our \cref{thm:uniquelyergodic} is hence an immediate corollary of the general
statement proved by C.~Fougeron in~\cite{Fougeron:Simplicial} for a broad class
of MCF.

\begin{remark}
	A result on genericity of the unique ergodicity property for
	Bruin-Troubetzkoy ITMs of infinite type was established
	in~\cite{BruinRadinger}. Namely, using a symbolic presentation of the ITMs,
	they established rather simple combinatorial conditions under which unique
	ergodicity holds for almost all (with respect to any renormalization
	invariant measure whose support set is $\gasket$) ITMs of infinite type.
	Both our proof and the one in~\cite{BruinRadinger} are based on the idea
	originated by Veech for IETs in~\cite{Veech:ue}. The key feature of our
	statement is the explicit construction of the ergodic measure $\mu$.
\end{remark}

Our approach also allows us to get another improvement for the result by Bruin
and Troubetzkoy. Namely, we prove the following estimations on the Hausdorff
dimension of the set $\gasket$ mentioned above: 

\begin{theorem}\label{thm:hdimboundaries} Let $\gasket$ be the set of parameters
	$(\alpha, \beta)$ yielding infinite type Bruin-Troubetzkoy ITMs. Then, its
	Hausdorff dimension can be bounded by
	\[
		1.5 \le \dim_H (\gasket) < 2.
	\]
	Moreover, the Hausdorff dimension of the set $\gasket$ is equal to its
	affinity dimension.
\end{theorem}

We refer to \cref{sec:HausdorffEstimate} for the definition of affinity
dimension. In the previous result, the upper bound follows from the application
of Fougeron's criterion, in a fashion similar to how we obtain
\cref{thm:uniquelyergodic}. More precisely, the upper bound follows from
Theorem~1.5 in~\cite{Fougeron:Simplicial}, which we can apply as our induction
satisfies the quickly escaping property, see \cref{lemma:quickscape} and
\cref{def:deg}. The lower bound is achieved by applying the strategy developed
in~\cite{JiaoLiPanXu:Applications} for another fractal, of rather similar
origin, called the \emph{Rauzy gasket}. The Rauzy gasket was widely studied in
the literature for several reasons, including symbolic dynamics, Arnoux-Rauzy
interval exchange transformations, pseudogroups of rotations and
$\mathbb{R}$-trees as well as Novikov's problem of asymptotic behavior of plane
sections of triply periodic surfaces, see \cite{DHS} for more details and
references. In fact, to prove that the Hausdorff dimension is equal to the
affinity dimension, we follow a very recent paper by N.~Jurga who obtains a
similar result for the Rauzy gasket~\cite{Jurga:Gasket}. 

For technical reasons, it will be convenient to change from the parameters
$(\alpha, \beta)$ to parameters $(a, b, c)$ in the $2$-dimensional symplex,
see~\eqref{eq:changeofcoords} on \cpageref{eq:changeofcoords}. This is
legitimate, since the change of coordinates is given by a linear and invertible
map, which preserves the Hausdorff dimension, since it is bi-Lipschitz. By an
abuse of notation, we will denote both sets by $\gasket$.

\begin{remark}
	The upper bound can be slightly improved using the strategy developed by
	Policott and Sewell for the Rauzy gasket~\cite{PollicottSewell},
	see~\cite{Zernikov:Diplom}. 
\end{remark}

\begin{figure}[t]
	\centering
	\includegraphics[width=.9\textwidth]{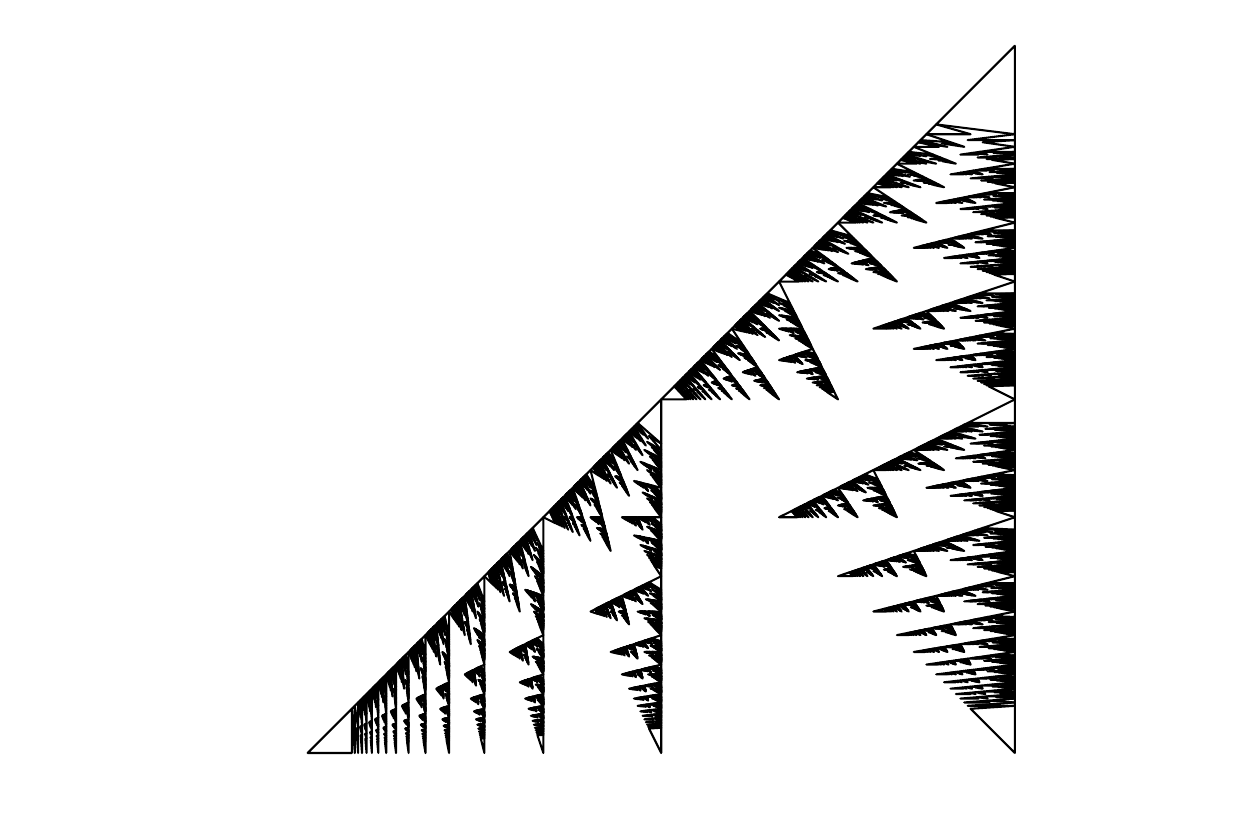}
	\caption{The Bruin-Troubetzkoy gasket.}
	\label{fig:BTGasket}
\end{figure}

In view of~\cref{thm:hdimboundaries}, we call the set $\gasket$ the
\emph{Bruin-Troubetzkoy gasket}, see \cref{fig:BTGasket,fig:BTGasket_simplex} on
\cpageref{fig:BTGasket_simplex} for a portrait of this set using the parameters
$(\alpha, \beta)$ and $(a, b, c)$, respectively. The topological similarity
between this fractal and the Rauzy gasket follows from the structure of the
renormalization algorithm we construct, which is quite similar to well-known and
well studied MCF algorithms, such as the Arnoux-Rauzy map, the Cassaigne
algorithm and the Arnoux-Rauzy-Poincar\'e algorithms, see~\cite{CLL} and the
references therein for the details. In order to reflect these features we give
to our algorithm a special name, we call it the \emph{Arnoux-Rauzy-Cassaigne
algorithm} (or ARC, for short). 

Using the MCF point of view, it is natural to compare the ergodic and spectral
properties of our algorithm with the ones of the above mentioned algorithms, and
the renormalization algorithm itself is the second dynamical system we are
looking at in the current paper. Our main result about the ARC MCF algorithm is
the following:

\begin{theorem}\label{thm:pisot}
	The cocycle defined by the ARC algorithm has almost always \emph{Pisot}
	Lyapunov spectrum. 
\end{theorem} 

To obtain it, we follow the ideas from~\cite{CLL}. Contrary to this result, it
is known that self-similar Bruin-Troubetzkoy ITMs of infinite type are very
often weakly mixing \cite{BruinRadinger, Mercat}. Moreover, we believe that weak
mixing for Bruin-Troubetzkoy ITMs is the typical behavior. In fact, we
conjecture that: 

\begin{conjecture}\label{weakmixingconjecture} 
	Almost all (with respect to the measure $\mu$ obtained in
	\cref{thm:uniquelyergodic}) Bruin-Troubetzkoy ITMs of infinite type are
	weakly mixing. 
\end{conjecture} 

Therefore, we have discovered an interesting phenomenon that did not appear
before in the situations that can be used as a natural references in our context
(generic IETs,  Arnoux-Rauzy IETs, and so on): we have a dynamical system which
is typically weakly mixing while the renormalization algorithm satisfies the
Pisot condition.  
This seems to be an interesting phenomenon, which warrants further investigation. 

\subsection*{Organization of the paper}
The paper is organized as follows: we start with the detailed description of the
renormalization algorithm (see \cref{sec:induction}); in the same section we
prove \cref{thm:uniquelyergodic} and the upper bound in
\cref{thm:hdimboundaries}. \cref{sec:Pisot} is devoted to the proof of
the \cref{thm:pisot}. Finally, the lower bound in the \cref{thm:hdimboundaries} is
proved in the \cref{sec:HausdorffEstimate}.

\subsection*{Acknowledgments}
We would like to thank Sebastien Labb\'e for interesting conversations and some
improvements for the first version of the paper; we are also grateful to Paul
Mercat for fruitful discussions. We thank Juan Galvis, Yessica Trujillo and Juan
Pablo Sierra for their help in preparing
\cref{fig:BTGasket,fig:BTGasket_simplex}. The first author thanks the Laboratory
of Cluster Geometry at HSE in Moscow for its hospitality. 

\section{Renormalization}\label{sec:induction} 

In this section we first describe the induction procedure for Bruin-Troubetzkoy
family of ITMs. Then, we apply it to get \cref{thm:uniquelyergodic}.

The induction we use is morally the same as the one used for interval exchange
transformations and for systems of isometries under the name \emph{Rauzy
induction}; as in these classical cases, the renormalization consists of
subsequent iterations of an algorithm, one step of which is the first return map
on the certain subinterval of the original support interval.

\subsection{Notation}
First, we change the notation in order to make the description of our family
more homogeneous. Namely, we introduce new parameters: if $\alpha>\beta$, we
have
\begin{equation}\label{eq:changeofcoords}
	\begin{split}
		a &= 1 - \alpha,\\
		b &=\alpha -\beta,\\
		c &= \beta.
	\end{split}
\end{equation}
Thus, $a + b + c = 1$ and 
\[
	\begin{split}
		T([0,a)) 	&= [1-a,1).\\
		T([a,a+b))  &= [1-b,1).\\
		T([a+b,1))  &= [0, c).
	\end{split}
\]
We always assume that $a$, $b$, and $c$ are \emph{rationally independent}, since
the subset of parameters where this condition does not hold comprises only
subsets of codimension 1 or 2. We remark that the parameters $a$, $b$, and $c$
correspond to the lengths of the intervals of continuity of the
Bruin-Troubetzkoy ITM $T_{\alpha,\beta}$.

Let us also enumerate the intervals of continuity of the map $T$ from the left
to the right; thus, the first interval is the one of the length $a$, the second
is the one of the length $b$ and the third is the one of the length $c$. The
vector that codes the order in which the intervals appear at the preimage of
$T$, is given by $\pi = (1,2,3).$ We denote the vector of lengths of the
intervals by $\vec{\lambda} = \transpose{(a,b,c)}$, where $\transpose{\cdot}$
denotes the transpose.

\subsection{The induction $\rauzy$}\label{sec:rauzy_induction}
To define our induction, we distinguish three cases. 

\begin{figure}[t]
	\centering
	\subfloat[][Case 1.\label{fig:Case1}]{\includegraphics[width=.4\textwidth]{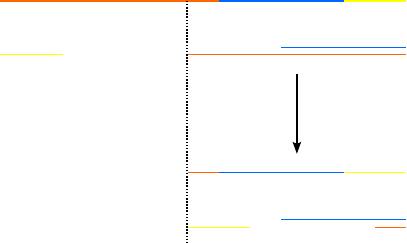}}
	\qquad\qquad
	\subfloat[][Case 3.\label{fig:Case3}]{\includegraphics[width=.4\textwidth]{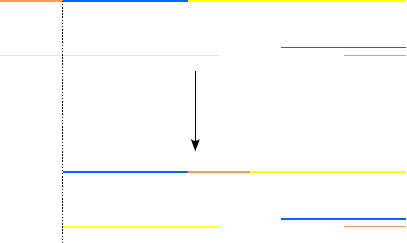}}
	\caption{The two cases of the $\rauzy$ induction not (immediately) yielding a finite type ITM.}
	\label{fig:RauzyInfiniteCases}
\end{figure}

{\bf Case 1: $a>b+c$.}\label{case1} We consider the first return map on the
subinterval $[b+c,1)$. It is an ITM in the same family with the following
lengths of intervals: 
\[
	\begin{split}
		a' &= a-b-c,\\
		b' &= b,\\
		c' &= c,
	\end{split}
\]
see \cref{fig:Case1}. We observe that the order of intervals does not change:
the interval of length $a'$ is still the leftmost, the interval of length $b'$
is in the middle, while the interval of length $c'$ is the rightmost. So, the
coding is again given by the vector $\pi' = \pi = (1,2,3)$.  

It is convenient to write that 
\[
	\vec\lambda = A\vec\lambda',
\]
where 
\[
	A = 
		\begin{pmatrix}
			1 & 1 & 1\\
			0 & 1 & 0\\
			0 & 0 & 1
		\end{pmatrix},
\]
and $\vec\lambda' = \transpose{(a',b',c')}$ is the vector of the new lengths of
the intervals, after the induction. We remark that the matrix $A$ of the
induction in this case coincides with the one defined by P.~Arnoux and G.~Rauzy
in~\cite{AR}. 

\begin{figure}[tb]
	\centering
	\includegraphics[width=.6\textwidth]{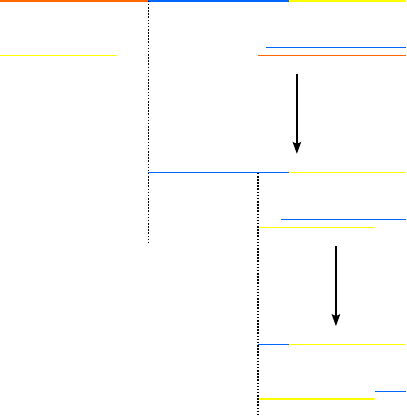}
	\caption{The case of the $\rauzy$ induction inducing a finite type ITM.}
	\label{fig:Case2}
\end{figure}

{\bf Case 2: $c<a<b+c$.}\label{case2} One can check that in this case the ITM
can be reduced to the ITM on 2 intervals and thus belongs to the finite case,
see \cref{fig:Case2}.

{\bf Case 3: $a<c$.}\label{case3} We consider the first return map to the
subinterval $[a,1)$. As result we get the following ITM:
\[
	\begin{split}
		T([a,a+b)) &= [1-b,1),\\
		T([a+b,2a+b)) &= [a-1,1),\\
		T([1-c+a,1)) &= [0,c-a).
	\end{split}
\] 
Then, the lengths change in the following way: 
\[
	\begin{split}
		a' &= a,\\
		b' &= b,\\
		c' &= c-a,
	\end{split}
\]
see \cref{fig:Case3}. As for Case 1, we use a matrix form: 
\[
	\vec\lambda = C_A\vec\lambda',
\]
where 
\[
	C_A = 
		\begin{pmatrix}
			1 & 0 & 0\\
			0 & 1 & 0\\
			1 & 0 & 1
		\end{pmatrix}.
\]

However, this case is very different from the first case, since the
\emph{order} of the intervals has changed: now the interval of length $b'$ is
the leftmost one, while the interval of length $a'$ is in the middle (the
third interval is always the rightmost). Note that the position of the intervals
in the image does not change: the image of the rightmost interval is always in
the left part and contains $0$, while the two other intervals are on the right
and contain the rightmost point of the support interval. So, the coding is given
by the vector $\pi' = (2,1,3)$. 

Now, since the induction consists in subsequent iterations of the
procedure described above, it is also necessary to analyze what happens when we
start with an ITM with coding $(2,1,3)$.

It is easy to see that if we start with an ITM with the combinatorial coding
given by the vector $\pi = (2,1,3)$ we get the symmetric picture. More
precisely, if the interval labelled by $2$ is longer than half of the support
interval, after the induction we still have the intervals in the order $\pi' =
(2,1,3)$. Whereas if the interval labelled by $3$ is longer than the rightmost
interval (in our case it is interval labelled by $1$), then the resulting ITM is
coded by $\pi' = (1,2,3)$ again. 

The matrices that are used to express the old lengths in terms of
new ones are slightly different: as an analogue of the matrix $A$ we get
\[
	B = 
		\begin{pmatrix}
			1 & 0 & 0\\
			1 & 1 & 1\\
			0 & 0 & 1
		\end{pmatrix}.
\]
and as an analogue of the matrix $C_A$ one has 
\[
	C_B = 
		\begin{pmatrix}
			1 & 0 & 0\\
			0 & 1 & 0 \\
			0 & 1 & 1
		\end{pmatrix}.
\]
Again, we remark that the matrix $B$ coincides with the one defined by P. Arnoux
and G. Rauzy in~\cite{AR}.

In all these situations we use general notation $\vec\lambda =
\mathcal{R}\vec\lambda'$.

\begin{figure}[bt]
	\begin{tikzpicture}[>=stealth]
		\node (left) [left=2cm]	{$(1, 2, 3)$};
		\node (right) [right=2cm] {$(2, 1, 3)$};
		\draw[->,>=latex] (left) to[bend left]
		node[above] {$C_A$}
		(right);
		\draw[->,>=latex] (right) to[bend left]
		node[above] {$C_B$}
		(left);
		\draw[->,>=latex] (right) to [out=30, in=330, looseness=4]
		node[right] {$B$}
		(right);
		\draw[->,>=latex] (left) to [out=150, in=210, looseness=4]
		node[left] {$A$}
		(left);
	\end{tikzpicture}
	\caption{The Rauzy graph $\rauzygraph$ of the induction $\rauzy$.}
	\label{fig:Rauzy_graph}
\end{figure}

Thus the induction process can be seen as a Markovian multidimensional fraction
algorithm. The diagram associated with this system is presented in
\cref{fig:Rauzy_graph}. In the Figure, we ignored Case 2 since it corresponds to
the finite case ITMs. We can think of Case 2 as a ``hole'', since if we enter in
it we escape from the world of infinite type ITMs.

The following lemma clearly holds: 

\begin{lemma}\label{infinitetype} 
The set $\gasket$ of rationally independent parameters that give rise to
infinite type Bruin-Troubetzoy ITMs coincides with the set of parameters for
which Case 2 of the induction $\rauzy$ never happens. That is, the parameters
which do not enter the hole during the induction procedure.
\end{lemma}

\subsection{The induction $\rauzy$ as a simplicial systems}\label{sec:simplicial}

In this and the next section, we make preliminary work and then apply a
criterion of Fougeron~\cite{Fougeron:Simplicial} to
obtain~\cref{thm:uniquelyergodic}. His criterion provides a general approach on
how to construct a measure on a fractal set if the set can be seen as an
attractor of an iterated function system coming from a Markovian
multidimensional continued fraction (MCF). The measure we construct can be seen
as a projection of the measure of maximal entropy on the canonical suspension
for the considered algorithm. 

As a by-product of Fougeron's approach we get an upper bound on the Hausdorff
dimension of the fractal $\gasket$ set as well as the unique ergodicity
statement for BT ITMs of infinite type\footnote{We will sometimes use the
abbreviation BT for Bruin-Troubetzkoy.}. The key notion for Fougeron's formalism
is the one of \emph{simplicial system}: a special graph associated with the
considered MCF. In this section we explain how to associate a simplicial system
with the induction algorithm. In the next one we check that all the conditions
of Fougeron's criterion are satisfied in our case. We refer
to~\cite{Fougeron:Simplicial} for more details and examples.

Let $G=(V,E)$ a graph labeled on an alphabet $\cA$ by a function $l\colon
E\to\cA$. We require that, for all $v\in V$, the restriction of $l$ to the
edges starting at $v$ in injective. Let $\|\cdot\|_1$ be the $\ell^1$-norm on
$\RR_+^\cA$ given by
$\|\lambda\|_1=\sum_{\alpha\in\cA}\lambda_\alpha$. Consider the simplex of
dimension $|\cA|-1$ defined by $\Delta=\{\lambda\in\RR_+^\cA :
\|\lambda\|_1=1\}$.

Given a vertex $v\in V$ let $v_{\text{out}}$ be the set of all edges going out
of $v$ in $G$. We define a partition of $\Delta$ whose atoms, for an edge $e\in
E$, are given by
\[
	\Delta^e=\{ \lambda\in\Delta : \lambda_{l(e)} < \lambda_\alpha, \text{ for all } \alpha\in l(v_{\text{out}}) \text{ and } \alpha\neq l(e) \}.
\]
To the edge $e$ we associate the matrix
\[
	M_e = \operatorname{Id} + \sum_{\substack{\alpha\in l(v_{\text{out}}) \\ \alpha\neq e}} E_{\alpha,l(e)},
\]
where $E_{a,b}$ is the elementary matrix with coefficient $1$ at the row $a$ and
column $b$ and all the other coefficients equal to $0$.

Define $\Delta^G=V\times\Delta$ and define the map $T\colon\Delta^G\to\Delta^G$
as follows. For all $\lambda\in\Delta^e$ with the edge $e$ starting at $v$ and
ending at $v'$, we say
\[
	T(v,\lambda)=(v', T_e(\lambda)),
\]
with
\[
	\begin{split}
		T_e\colon\Delta^e &\to\Delta\\
			\lambda &\mapsto \frac{M_e^{-1}\lambda}{\|M_e^{-1}\lambda\|}.
	\end{split}
\]
We say that $G$ is a \emph{simplicial system} and the map $T$ is its associated
\emph{win-lose induction} (see below for an explanation of this name).

From these definitions it is clear that the (projectivization of) the induction
$\rauzy$ we described in \cref{sec:rauzy_induction} is obtained from the
win-lose induction associated to the simplicial system defined by the graph $G$
in \cref{fig:SimplicialRauzy} by considering the first return map to the white
vertices.

\begin{figure}[tb]
	\center
	\begin{tikzpicture}[shorten >=1pt,node distance=2cm,on grid,auto]
	   \node (11) {$\circ$};
	   \node (12) [right = of 11] {};
	   \node (13) [right = of 12] {$\circ$};
	   \node (21) [below = of 11] {$\bullet$};
	   \node (22) [below = of 12] {$\bullet$};
	   \node (23) [below = of 13] {$\bullet$};

	   \path[->] (11) edge [bend left=10] node {$1$} (13);
	   \path[->] (13) edge [bend left=10] node {$2$} (11);
	   \path[->] (11) edge [bend left=20] node {$3$} (21);
	   \path[->] (13) edge [bend left=20] node {$3$} (23);
	   \path[->] (21) edge [bend left=20] node {$2$} (11);
	   \path[->] (23) edge [bend left=20] node {$1$} (13);
	   \path[->] (21) edge 				  node [below] {$1$} (22);
	   \path[->] (23) edge 				  node  {$2$} (22);
	
	\end{tikzpicture}
	\caption{The induction $\rauzy$ as a simplicial system.}
	\label{fig:SimplicialRauzy}
\end{figure}

\subsection{Proof of \cref{thm:uniquelyergodic}}  
We will use the machinery developed in~\cite{Fougeron:Simplicial}, and
especially the result in Section~4 of that paper, to prove our first main
result, namely \cref{thm:uniquelyergodic}: the existence of an
$\rauzy$-invariant measure supported on $\gasket$, with respect to which almost
every Bruin-Troubetzkoy ITM is uniquely ergodic. Hence, we need to check that
after a sufficiently large amount of steps we get a positive matrix for the
induction. So we are interested in what happens when we apply several
consecutive steps of the algorithm described above. We recall the terminology
used for the classical Rauzy induction for IETs. In each iteration of the
induction, the longest interval, i.e., the one that gets cut, will be called the
\emph{winner} and the shortest one is called the \emph{loser}. As with IETs, we
name the intervals using the corresponding letter. With this convention, in Case
1 the $a$-interval is the winner and the $c$-interval is the loser, whereas in
Case 3 it is the opposite. We will sometimes simply say that the letter (and not
the corresponding interval) is the winner or the loser. The following statement
follows from \cref{infinitetype}: 

\begin{lemma}\label{winlose} 
	A Bruin-Troubetzkoy ITM is of infinite type if and only if each letter wins
	and loses infinite number of times. 
\end{lemma} 

A corollary of this result is that the simplicial system associated to the
induction $\rauzy$ is \emph{quickly escaping} in the sense
of~\cite{Fougeron:Simplicial}. Actually, we will show a stronger, yet easier to
check property: strong non-degeneracy. Let us recall this terminology.

Given a simplicial system $G$ and $\cL\subset\cA$, we denote by $G_\cL$
the subgraph of $G$ with the same set of vertices $V$ and whose edges are as
follows. For a vertex $v\in V$:
\begin{itemize}
	\item if at least one outgoing edge from $v$ is labeled in $\cL$ then the
	set of outgoing edges in $G_\cL$ is
	\[
		v_\text{out}^\cL = \{e\in v_{\text{out}} : l(e)\in\cL\};
	\]
	\item otherwise $v_{\text{out}}^\cL=v_{\text{out}}$.
\end{itemize}

\begin{definition}[Strongly non-degenerating]
	\label{def:deg}
	We say that a simplicial system is \emph{strongly non-degenerating} if
	\begin{enumerate}
		\item for all vertices every letter wins and loses in almost every path
		with respect to Lebesgue measure,
		\item for all $\emptyset\subsetneq\cL\subsetneq\cA$ and all vertex
		$v$ in a strongly connected components $\mathscr{C}$ of $G_\cL$ then
		either
			\begin{itemize}
				\item $|l(v_\text{out})\cap\cL|\leq 1$;
				\item or there is a path from $v$ in $G$, labeled in $\cL$,
				which leaves $\mathscr{C}$.
			\end{itemize}
	\end{enumerate}
\end{definition}

Then, we have the following result.

\begin{lemma}\label{lemma:quickscape}
	The simplicial system associated to the induction $\rauzy$ is quickly
	escaping.
\end{lemma}

\begin{proof}
	Theorem~3.13 of~\cite{Fougeron:Simplicial} states that a strongly
	non-degenerating simplicial system is quickly escaping. By \cref{winlose},
	the first condition of \cref{def:deg} holds. The second condition can be
	easily checked: for every choice of $\cL$ we obtain a graph $G_\cL$ which
	has only one strongly connected component in which every vertex $v$
	satisfies $|l(v_\text{out})\cap\cL|\leq 1$.
\end{proof}

Thus, given that $\lambda = (a,b,c)$ is the vector of lengths of the intervals
and $\lambda^{(n)} = (a^{(n)}, b^{(n)}, c^{(n)})$ is the vector of lengths of
the intervals after the application of $n$ steps of the induction, we have
$\lambda = \mathcal{R}^{(n)} \lambda',$ with $\mathcal{R}^{(n)}
=\mathcal{R}(k_1,k_2,k_3,\cdots) = A^{k_1}C_{A}B^{k_2}C_{B}A^{k_3}\cdots,$ where
$k_1, k_2, \cdots \in\mathbb{N}, \sum_i{k_i} = n$. The same formula holds if we
let $n$ go to infinity, and consider infinite paths on the graph in
\cref{fig:Rauzy_graph}.

We stress that the coefficients $k_i$ can be equal to $0$. However, since
applying the matrix $A$ implies cutting of $a$-intervals, we must have that, for
an infinite number of steps, $k_{2i+1} > 0$. Similarly, since the application of
the matrix $B$ implies cutting the $b$-interval, we must have that, for an
infinite number of steps, $k_{2i} > 0$. Hence, an ITM is of infinite type if and
only if we have infinitely often that even and odd $k_i$'s are \emph{strictly}
positive.

Now one can check that any $\rauzy$ that contains $A$ and $B$ in positive
powers together with $C_A$ and $C_B$ has strictly positive entries. Therefore,
the following lemma holds: 

\begin{lemma}
There exists a special acceleration of the induction described above. 
\end{lemma} 

The definition of \emph{special acceleration} can be seen
in~\cite[Remark~1]{FS}. Morally, it means a first return map to some subsimplex
compactly contained in the parameter space. Exploiting the machinery of
simplicial systems introduced in~\cite{Fougeron:Simplicial} and related results
from~\cite{FS}, we easily obtain \cref{thm:uniquelyergodic}.

\begin{proof}[Proof of \cref{thm:uniquelyergodic}]
	The simplicial system associate to $\rauzy$ is \emph{uniformly expanding}
	by~\cite[Proposition~4.1]{Fougeron:Simplicial} and therefore ergodic thanks
	to~\cite[Corollary~4.4]{Fougeron:Simplicial}. \cref{lemma:quickscape}
	guarantees that the simplicial system is quickly escaping and thus, by
	Theorem~4.24 in \cite{Fougeron:Simplicial}, we obtain the natural measure
	$\mu$ that induces the measure of maximal entropy on the natural suspension.
	Therefore the set of parameters which follow the same path (generic for
	$\mu$) is a single point, and so by the standard argument of Veech
	(see~\cite{Veech:ue}) we conclude that the original ITM is uniquely ergodic.
	This completes the proof of \cref{thm:uniquelyergodic}.
\end{proof}

As a corollary of the previous result, we obtain an upper estimate on the
Hausdorff dimension of the parameters yielding Bruin-Troubetzkoy ITMs of
infinite type.

\begin{corollary}
	The set $\gasket$ of parameters that give rise to the infinite type
	Bruin-Troubetzkoy ITMs has Hausdorff dimension strictly smaller than $2$.
\end{corollary}

\begin{proof}
    At any vertex of the graph underlying the simplicial system associated to
    the $\rauzy$ induction, the set of lengths vectors for which the induction
    does not stop after finitely many iterations forms a subset of the simplex
    (which is denoted $\Delta(F)$ in~\cite{Fougeron:Simplicial}). By
    \cref{lemma:quickscape}, the simplicial system is quickly escaping and so we
    can use~\cite[Theorem~1.5]{Fougeron:Simplicial} to conclude.
\end{proof}

\begin{figure}[t]
	\centering
	\includegraphics[width=.9\textwidth]{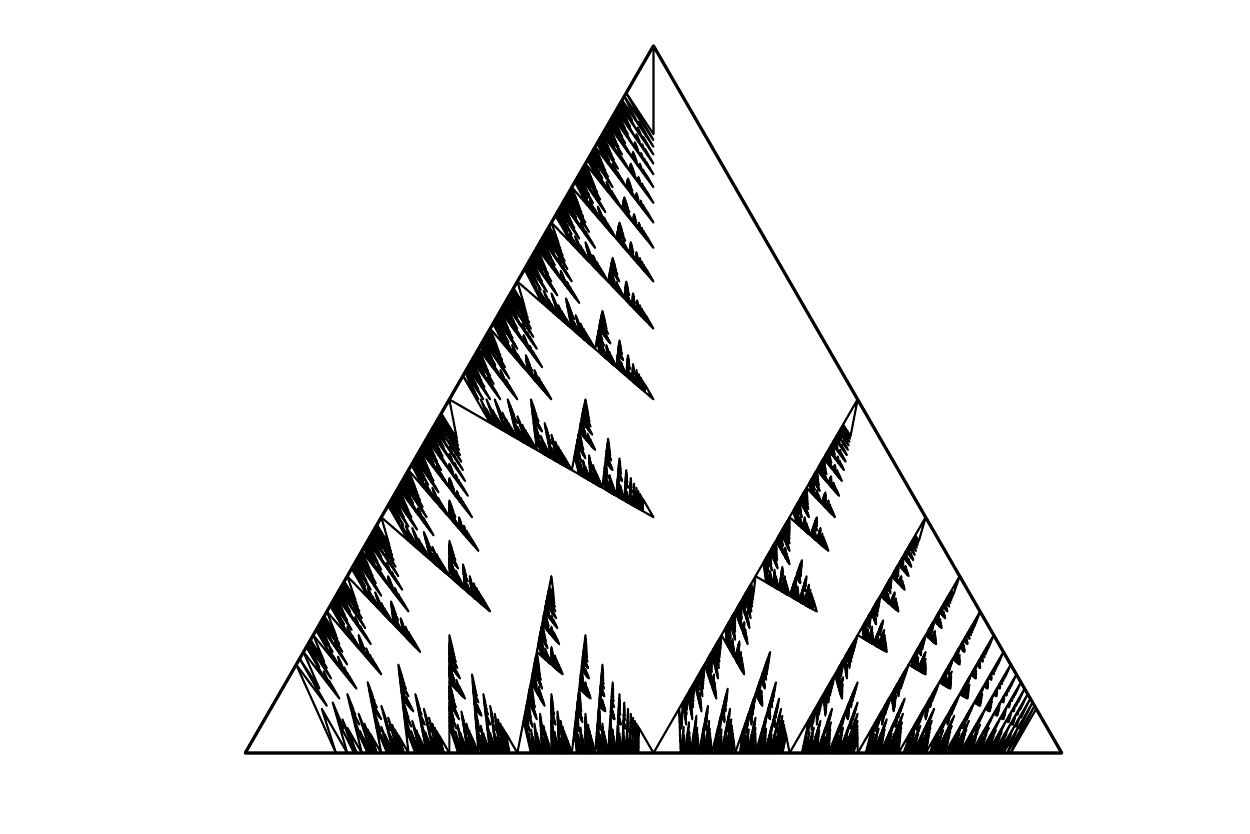}
	\caption{The Bruin-Troubetzkoy gasket using the simplicial coordinates $(a, b, c)$.}
	\label{fig:BTGasket_simplex}
\end{figure}

We will obtain a lower bound, using thermodynamical formalism, in
\cref{sec:HausdorffEstimate}. In \cref{fig:BTGasket_simplex}, we represent the
Bruin-Troubetzkoy gasket using the parameters $(a, b, c)$ instead of $(\alpha,
\beta)$ as in \cref{fig:BTGasket}. 

\subsection{Recovering Bruin and Troubetzkoy's Gauss map} 
We now show that our induction can be accelerated to recover the Gauss map of
Bruin and Troubetzkoy.

\begin{proposition}\label{prop:RauzyGauss}
	There exists an acceleration of $\rauzy$ such that, after rescaling the
	original interval, the induced transformation is the one obtained via the
	Gauss map of Bruin and Troubetzkoy.
\end{proposition}

\begin{proof}
	We begin by recalling the definition of the Gauss map. If $T_{\alpha,
	\beta}$ is a Bruin-Troubetzkoy ITM, then we define the ITM $T_{\alpha',
	\beta'}$ where
	\begin{equation}\label{eq:GaussMap}
		(\alpha', \beta') = \biggl(\frac{\beta}{\alpha}, 
			\frac{\beta-1}{\alpha} + \biggl\lfloor \frac{1}{\alpha} \biggr\rfloor\biggr),
	\end{equation}
	where $\lfloor \cdot \rfloor$ is the (lower) integer part. Let us recall
	that we can recover the parameter $\alpha$ and $\beta$ from the length ones
	using that $\alpha = b + c$ and $\beta = c$, if the intervals are in the
	order $(1,2,3)$, and similarly, replacing $b$ by $a$ in the other case.
	
	We observe that, by definition of the $\rauzy$  induction, the first case is
	repeated $n$ times, with $n \ge 0$ given by
	\[
		n 	= \biggl\lfloor \frac{a}{b+c} \biggr\rfloor
			= \biggl\lfloor \frac{1-\alpha}{\alpha} \biggr\rfloor 
			= \biggl\lfloor \frac{1}{\alpha} \biggr\rfloor  - 1.
	\]
	Then, we are in the third case, and we change the order of the intervals.

	After the above steps, the three intervals of continuity are of lengths
	\[
		\begin{split}
			a' &= a - n(b+c),\\
			b' &= b,\\
			c' &= c - (a- n(b+c)). 
		\end{split}
	\]
	The total length is $a'+b'+c' = b+c = \alpha$. So, if we renormalize,
	dividing the interval by $\alpha$, rescaling it to be of length $1$, we see
	that
	\[
		\alpha' = \frac{a'+c'}{\alpha} = \frac{c}{\alpha} = \frac{\beta}{\alpha}
	\]
	Moreover, since
	\[
		c' = \beta - 1 +\alpha + \biggl( \biggl\lfloor \frac{1}{\alpha} \biggr\rfloor - 1\biggr) \alpha = \beta - 1 + \biggl\lfloor \frac{1}{\alpha} \biggr\rfloor \alpha,
	\]
	we have that
	\[
		\beta' = \frac{c'}{\alpha} = \frac{\beta - 1}{\alpha} + \biggl\lfloor \frac{1}{\alpha} \biggr\rfloor.
	\]
	The above formulas agree with~\eqref{eq:GaussMap} and so we are done.
\end{proof}

\section{Pisot property for the ARC algorithm}\label{sec:Pisot}
\subsection{The ARC multidimensional continued fraction
algorithm}\label{sec:ARC_MCF} 

The induction $\rauzy$ introduced in the previous section defines a
multidimensional continued fraction algorithm (or MCF algorithm, for short),
which we call the \emph{Arnoux-Rauzy-Cassaigne} (or ARC) MCF algorithm, as we
will now explain. We can naturally act by the matrices of the $\rauzy$ induction
to the standard $2$-dimensional simplex $\Delta = \Delta^2 = \{(x,y,z): x,y,z\ge
0, x+y+z=1\}$. In formulas, we have:
\begin{equation}\label{eq:ARC_IFS}	
	\begin{split}
		f_A (x,y,z) &= \biggl(\frac{1}{2-x}, \frac{y}{2-x}, \frac{z}{2-x}\biggr),\\
		f_B (x,y,z) &= \biggl(\frac{x}{2-y}, \frac{1}{2-y}, \frac{z}{2-y}\biggr),\\
		f_{C_A} (x,y,z) &= \biggl(\frac{1-y}{2-x-y}, \frac{y}{2-x-y}, \frac{z}{2-x-y}\biggr),\\
		f_{C_B} (x,y,z) &= \biggl(\frac{x}{2-x-y}, \frac{1-x}{2-x-y}, \frac{z}{2-x-y}\biggr).
	\end{split}	
\end{equation}

We will now show that the cocycle defined by the ARC MCF algorithm has negative
second Lyapunov exponent (see below for the relevant definitions),
following~\cite{CLL}; in the terminology of~\cite{Lagarias}, the MCF algorithm
is \emph{strongly convergent}. Then, the Lyapunov exponents satisfy $\lambda_1 >
0 > \lambda_2$, which is the Pisot condition, see, e.g.,
\cite[Definition~2.1]{BST}. 

We remark that the measure $\mu$ obtained in the proof of
\cref{thm:uniquelyergodic} induces a measure, which, slightly abusing the
notation, we still denote $\mu$, on set of walks along the graph $\rauzygraph$.
We will call $1$ the state corresponding to the permutation $(1, 2, 3)$ and $2$
the permutation $(2, 1, 3)$. Thus, for instance, the path $112$ corresponds to
the application of the matrices $AC_A$ of the induction. As usual, let $[a_1,
a_2, \dotsc, a_n]$ be the cylinder formed by all the words in $\{1, 2\}^\NN$
that begin with the letters $a_1a_2\cdots a_n$. 

We now recall the definition of Lyapunov exponents in the present context (along
with the relevant notation). Given an infinite word $w\in\{1,2\}^\NN$, we
consider the product of matrices corresponding to the path on the graph
$\rauzygraph$ described by $w$:
\[
	X_n(w) = X_{w_0w_1} X_{w_1w_2} \cdots X_{w_{n-2}w_{n-1}},
\]
for $n\ge 1$. This forms a cocycle, as
\[
	X_{m+n}(w) = X_m(w) X_n(\sigma^m w),
\]
where $\sigma\colon\{1,2\}^\NN\to\{1,2\}^\NN$ is the shift map. Since all the
matrices of the $\rauzy$ induction are invertible, this cocycle is
log-integrable with respect to the measure $\mu$ constructed in
\cref{thm:uniquelyergodic}:
\[
	\int_{\{1,2\}^\NN} \log \max\bigl\{\|X_1(w)\|, \|X_1(w)^{-1}\|\bigr\} \, d\mu(w) < \infty.
\]
Hence there are ($\mu$-almost everywhere) well defined Lyapunov exponents:
\[
	\lambda_1 = \lim_{n\to\infty} \frac{\log{\|X_n(w)\|}}{n},
	\qquad
	\lambda_1 + \lambda_2 = \lim_{n\to\infty} \frac{\log{\|\wedge^2 X_n(w)\|}}{n},
	\qquad
	\lambda_3 = -(\lambda_1 + \lambda_2)
\]
where $\wedge$ is, as usual, the exterior product of matrices, and the last
equality follows from the fact that the matrices of the induction have
determinant $1$. Moreover, since, by unique ergodicity, the nested cones $X_1
X_2 \cdots X_n = X_{[1,n)}\RR^3_{\ge 0}$, where we have suppressed the
dependence on $w$, converge for $\mu$-almost all words $w$, to a line $\RR_{0}
f$, for some $f\in\RR^3$, we have the following characterization of the second
Lyapunov exponent will be useful later on:
\[
	\lambda_2 = \lim_{n\to\infty} \frac{\log{\Bigl\|\transpose{X_n(w)}|_{f^\perp}\Bigr\|}}{n},
\]
where $\transpose{X}$ denotes the transpose matrix and $f^\perp$ denotes the
orthogonal complement, with respect to usual inner product on $\RR^3$, of the
vector $f$.

\subsection{Some technical lemmas}
We begin with a general lemma about matrices and norms, whose proof can be found
in~\cite[Lemma~3.5]{CLL}.

\begin{lemma}\label{lemma:submultiplicative}

	Let $X$ and $Y$ non negative $d\times d$ real matrices, such that $XY\neq
	0$. Let $\|\cdot\|$ denote any seminorm on $\RR^d$ which is a norm on every
	$f^\perp$, with $f\in (X\RR^d_{\ge 0} \cup Y\RR^d_{\ge 0}) \setminus \{0\}$.
	We have
	\[
		\Bigl\| \transpose{XY} \Bigr\|^{XY\RR^d_{\ge 0}} \le 
		\Bigl\| \transpose{X} \Bigr\|^{X\RR^d_{\ge 0}} \Bigl\| \transpose{Y} \Bigr\|^{Y\RR^d_{\ge 0}}.
	\]
\end{lemma}

From the definitions, we obtain:

\begin{lemma}\label{lemma:products_of_M}
	Let $(M_n)_{n \in \NN} \in \{A, B, C_A, C_B\}^\NN$ be a sequence of
	matrices. If $(M_n)$ contains infinitely many occurrences of $A$, $B$, $C_A$
	and $C_B$, then there exists an increasing sequence of integers $(n_m)_{m
	\in \NN}$, such that $n_0 = 0$ and
	\[
		M_{[n_m, n_{m+1})} \in \{A^k C_A, B^k C_B: k\in\NN\}.
	\]
\end{lemma}

Following~\cite{CLL}, we introduce the seminorm $\|\cdot\|_D$ on $\RR^3$ given
by
\[
	\|v\|_D = \max v -\min v = \max_{i=1,2,3} v_i - \min_{i=1,2,3} v_i.
\]
We remark that the seminorm is invariant under addition of constant vectors.
Moreover, for any $f\in\RR^3_{\ge 0}$, the restriction of the seminorm to
$f^\perp$ is a genuine norm, see~\cite[Lemma~3.7]{CLL}. The same result implies
that, when restricted to $f^\perp$, for a given $f\in\RR^3_{\ge 0}$, we obtain a
norm on $3\times 3$ matrices, which is invariant under addition of constant
vectors, and comparable with the infinity norm. More precisely,
in~\cite[Lemma~3.7]{CLL} it is proven that: 
\begin{equation}\label{eq:Dnorm_vs_infinitynorm}
	\frac{1}{2}\bigl\| M|_{f^\perp} \bigr\|_D \le \bigl\| M|_{f^\perp} \bigr\|_\infty 
		\le 2 \, \bigl\| M|_{f^\perp} \bigr\|_D
\end{equation}

With this notation, we have

\begin{lemma}[\protect{\cite[Lemma~3.8]{CLL}}]\label{lemma:hyperplanes}
	Let $M$ be a $3\times 3$ positive and invertible real matrix. Consider the
	set $\cH$ of hyperplanes orthogonal to some vector in
	\[
		S = (M\cE \cup (\cE - \cE) \cup M(\cE - \cE)) \setminus \{0\},
	\]
	where $\cE = \{e_1, e_2, e_3\}$ and $\cD$ be the finite union of
	one-dimensional intersections of two hyperplanes of $\cH$:
	\[
		\cD = \bigcup_{h_1, h_2 \in \cH} h_1 \cap h_2.
	\]
	Then, the maximal value of the norm $\|\cdot\|_D$, for all restrictions
	orthogonal to a positive vector in the cone $M\RR^3_{\ge 0}$, is attained at
	some vector $v$ in $\Bigl(\cD \setminus \pm \transpose{M^{-1}} \RR^3_{>0}\Bigr)
	\setminus \{0\}$, i.e.,
	\[
		\Bigl\|\transpose{M}\Bigr\|_D^{M\RR^3_{\ge 0}} 
		:= \sup_{f\in M\RR^3_{\ge 0} \setminus \{0\}} \Bigl\|\transpose{M}|_{f^\perp}\Bigr\|_D
		 = \max_{v\in \bigl(\cD \setminus \pm \transpose{M^{-1}} \RR^3_{>0}\bigr) \setminus \{0\}} 
		 	\frac{\Bigl\|\transpose{M}v\Bigr\|_D}{\|v\|_D}.
	\]
\end{lemma}

We will use the next lemma to control the norm introduced above for specific
products of the matrices of our MCF.

\begin{lemma}\label{lemma:norm_one}
	For every $k\ge 0$ we have that 
	\[
		\Bigl\|\transpose{(A^kC_A)}\Bigr\|_D^{\RR_{\ge 0}^3} = 1, \qquad
		\Bigl\|\transpose{(B^kC_B)}\Bigr\|_D^{\RR_{\ge 0}^3} = 1.
	\]
\end{lemma}

\begin{proof}
	We will prove the lemma only for $A^kC_A$, the other case being similar.
	Choose a vector $f \in \RR^3_{\ge 0} \setminus \{0\}$ and let $v = (v_1,
	v_2, v_3) \in f^\perp$. Using the invariance of the norm $\|\cdot\|_D$ under
	addition by constant vectors, we obtain, for every $k\ge 0$
	\[
		\Bigl\|\transpose{(A^kC_A)}\Bigr\|_D^{\RR_{\ge 0}^3} = 
		\left\| 
			\begin{pmatrix}
				k+1 & 0 & 1 \\
				k 	& 1 & 0 \\
				k	& 0 & 1
			\end{pmatrix}
			\begin{pmatrix}
				v_1 \\ 
				v_2 \\
				v_3
			\end{pmatrix}
		\right\|_D
		=
		\left\|
			\begin{pmatrix}
				1 & 0 & 0 	\\
				0 & 1 & -1 	\\
				0 & 0 & 0
			\end{pmatrix}
			\begin{pmatrix}
				v_1 \\
				v_2 \\
				v_3
			\end{pmatrix}
		\right\|_D.
	\]

	Let $v' = (v_1, v_2 - v_3, 0)$, then 
	\[
			\min v \le \min v' \le 0 \le \max v' \le \max v.
	\]
	We remark that, since we work in $3$ dimensions, we have that $\|(v_1, v_2,
	v_3)\|_D = \max\{|v_1 - v_2|, |v_1 - v_3|, |v_2 - v_3|\}$, so the previous
	equation implies that $\|v'\|_D \le \|v\|_D$. Hence
	\[
			\frac{\Bigl\|\transpose{(A^k C_A)}v\Bigr\|_D}{\|v\|_D} \le 1,
	\] 
	for every $v\in f^\perp$ and $k\ge 0$. To obtain the equality, we observe
	that, since $f$ is non zero and non negative, we can take a vector $v = (a,
	-b, 0)$, for $a$, $b>0$ inside $f^\perp$. For this vector, the direct
	computation yields
	\[
		\Bigl\|\transpose{(A^k C_A)}(a, -b, 0)\Bigr\|_D = \|(a, -b, 0)\|_D,
	\]
	as we wanted.
\end{proof}

We need to take care of a ``base case'' before we can do the general one.

\begin{lemma}\label{lemma:AC_ABC_B}
	Let $M = AC_A BC_B$ or $M = BC_B AC_A$. Then
	\[
		\Bigl\| \transpose{M} \Bigr\|_D^{M\RR^3_{\ge 0}} \le \frac{4}{5}.
	\]
\end{lemma}

\begin{proof}
	We will prove the result in the case $M = AC_A BC_B$, the other case is
	symmetric. By direct computation:
	\[
		M = AC_ABC_B =
		\begin{pmatrix}
			3 & 3 & 2 \\
			1 & 2 & 1 \\
			1 & 1 & 1
		\end{pmatrix},
		\qquad
		M^{-1} =
		\begin{pmatrix}
			1  & -1 & -1 \\
			0  & 1 & -1 \\
			-1 & 0  & 3
		\end{pmatrix}.
	\]
	Given $z = (a, b, c)$ we can compute
	\[
		\Bigl\|\transpose{M}z\Bigr\|_D = \Bigl\|\transpose{(3a+b+c, 3a+2b+c, 2a+b+c)}\Bigr\|_D 
			= \Bigl\|\transpose{(0, b, a+b+c)}\Bigr\|_D.
	\]
	We now construct the set $\cH$ as in \cref{lemma:hyperplanes}. By a direct
	computation:
	\begin{align*}
		Me_1 &= (3, 1, 1), & e_1 - e_3 &= (1, 0, -1), & M(e_1 - e_3) &= (1, 0, 0), \\
		Me_2 &= (3, 2, 1), & e_1 - e_2 &= (1, -1, 0), & M(e_1 - e_2) &= (0, -1, 0), \\
		Me_3 &= (2, 1, 1), & e_2 - e_3 &= (0, 1, -1), & M(e_2 - e_3) &= (1, 1, 0).
	\end{align*}
	Hence, $\cH$ is made of nine hyperplanes. Again by \cref{lemma:hyperplanes},
	we need to consider vectors $z \in \Bigl(\cD \setminus \pm \transpose{M^{-1}}
	\RR^3_{>0}\Bigr) \setminus \{0\}$. The relevant computations are in
	\cref{table:CLL}, from which the result follows.
\end{proof}

\begin{table}[tp]
	\[
		\begin{array}{ccccccc}
			\toprule
			u & v & z = u \wedge v & \transpose{M}z &\|z\|_D & \|\transpose{M}z\|_D \\
			\midrule
			Me_1 & Me_2 		& (-1, 0, 3) 	& (0, 0, 1) 	& 4 & 1 \\
			Me_1 & Me_3 		& (0, -1, -1) 	& (0, -1, 0) 	& 2 & 1 \\
			Me_1 & e_1 - e_3 	& (-1, 4, -1) 	& (0, 4, 1) 	& 5 & 4 \\
			Me_1 & e_1 - e_2 	& (1, 1, -4) 	& (0, 1, -1) 	& 5 & 2 \\
			Me_1 & e_2 - e_3 	& (-2, 3, 3) 	& (0, 3, 2)  	& 5 & 3 \\
			Me_1 & M(e_1 - e_2) & (1, 0, -3) 	& (0, 0, -1) 	& 4 & 1 \\
			Me_1 & M(e_2 - e_3) & (-1, 1, 2) 	& (0, 1, 1) 	& 3 & 1 \\
			Me_1 & M(e_1 - e_3) & (0, 1, -1) 	& (0, 1, 0) 	& 2 & 1 \\
			Me_2 & Me_3 		& (1, -1, -1) 	& (1, 0, 0) 	& 2 & 1 \\
			Me_2 & e_1 - e_3 	& (-2, 4, -2) 	& (-4, 0, 2) 	& 6 & 4 \\
			Me_2 & e_1 - e_2 	& (1, 1, -5) 	& (-1, 0, -2) 	& 6 & 2 \\
			Me_2 & e_2 - e_3 	& (-3, 3, 3) 	& (-3, 0, 0) 	& 6 & 3 \\
			Me_2 & M(e_1 - e_2) & (1, 0, -3)	& (0, 0, -1) 	& 4 & 1 \\
			Me_2 & M(e_2 - e_3) & (-1, 1, 1) 	& (-1, 0, 0) 	& 2 & 1 \\
			Me_2 & M(e_1 - e_3) & (0, 1, -2) 	& (-1, 0, -1) 	& 3 & 1 \\
			Me_3 & e_1 - e_3 	& (-1, 3, -1) 	& (-1, 2, 0) 	& 4 & 3 \\
			Me_3 & e_1 - e_2 	& (1, 1, -3) 	& (1, 2, 0) 	& 4 & 2 \\
			Me_3 & e_2 - e_3 	& (-2, 2, 2) 	& (-2, 0, 0) 	& 4 & 2 \\
			Me_3 & M(e_1 - e_2)	& (1, 0, -2) 	& (1, 1, 0) 	& 3 & 1 \\
			Me_3 & M(e_2 - e_3) & (-1, 1, 1) 	& (-1, 0, 0) 	& 2 & 1 \\
			Me_3 & M(e_1 - e_3) & (0, 1, -1) 	& (0, 1, 0)  	& 2 & 1 \\
			\bottomrule
		\end{array}
	\]
	\caption{The computations involved in \cref{lemma:AC_ABC_B}. We only wrote
	the values of $u$ and $v$ which yield a $z = u \wedge v$ in
	$\RR^3\setminus\transpose{M}^{-1}\RR^3_{>0}$.}\label{table:CLL}
\end{table}

The main technical result is the following

\begin{lemma}\label{lemma:measure_cylinder}
	Let $\mu$ the measure on $\{1,2\}^\NN$ obtained in
	\cref{thm:uniquelyergodic}. For every $\epsilon > 0$, there exists an $N$
	such that, for every $n > N$ and $\mu$-almost every sequences $(M_n)_{n \in
	\NN} \in \{A, B, C_A, C_B\}^\NN$, we have
	\[
		\Bigl\| \transpose{M}_{[0,n)}|_{f^\perp}\Bigr\|_\infty 
		\le (n+1) \biggl(\frac{4}{5}\biggr)^{\frac{1}{8}n(\mu([112211221])-\epsilon)-\frac{1}{8}},
	\]
	where $\bigcap_{n\in\NN} M_{[0,n)}\RR^3_{\ge 0} = \RR_{\ge 0} f$.
\end{lemma}

\begin{proof}
	We begin by recalling that, by its construction, the measure $\mu$ assigns
	positive measure to every cylinder. In particular, $\mu([1])$, $\mu([2])>0$.
	By ergodicity of $\mu$, $\mu$-almost every sequence of matrices
	$(M_n)_{n\in\NN}$ contains infinitely often each one of the matrices $A$,
	$B$, $C_A$ and $C_B$. Then, by \cref{lemma:products_of_M}, there is an
	increasing sequence $(n_i)_{i\in\NN}$ such that $n_0 = 0$ and
	\[
		N_i = M_{[n_i, n_{i+1})} \in \{A^k C_A, B^k C_B: k\in\NN\}
	\]
	for all $i$. For all positive $n$, there exists a unique $m\in\NN$ such that
	$n_m \le n-1 < n_{m+1}$. Let $g = M_{[0,n)}^{-1} f$, then by
	\cref{lemma:submultiplicative}, we obtain
	\[
		\begin{split}
			\Bigl\| \transpose{M}_{[0,n)}|_{f^\perp}\Bigr\|_\infty &\le \Bigl\| \transpose{M}_{[0,n)}\Bigr\|_\infty^{M_{[0,n)} \RR^3_{\ge 0}} \\
				&\le \Bigl\| \transpose{M}_{[n_m,n)}\Bigr\|_\infty^{M_{[n_m, n)} \RR^3_{\ge 0}} \cdot \Bigl\| \transpose{M}_{[0,n_m)}\Bigr\|_\infty^{M_{[0, n_m)} \RR^3_{\ge 0}}\\
				&\le \Bigl\| \transpose{M}_{[n_m,n)}\Bigr\|_\infty \cdot \Bigl\| \transpose{M}_{[0,n_m)}\Bigr\|_\infty^{M_{[0, n_m)} \RR^3_{\ge 0}}
		\end{split}
	\]
	Since $M_{[n_m,n)}$ is of the form
	\[
		A^k =
		\begin{pmatrix}
			1 & k & k \\
			0 & 1 & 0 \\
			0 & 0 & 1
		\end{pmatrix},
		\qquad
		B^k =
		\begin{pmatrix}
			1 & 0 & 0 \\
			k & 1 & k \\
			0 & 0 & 1
		\end{pmatrix},
	\]
	for some $k\in\NN$, and $\Bigl\| \transpose{A^k}\Bigr\|_\infty = \Bigl\|
	\transpose{B^k}\Bigr\|_\infty = k+1$, we have
	\[
		\Bigl\| \transpose{M}_{[n_m,n)}|_{f^\perp}\Bigr\|_\infty \le n - n_m + 1 \le n + 1.
	\]

	Now, we deal with the second term: $\Bigl\|
	\transpose{M}_{[0,n_m)}\Bigr\|_\infty^{M_{[0, n_m)} \RR^3_{\ge
	0}}$, with $M_{[0,n_m)} = \prod_{i=0}^{m-1} N_i = N_{[0,m)}$. Let $J_m$ be
	the set of indices $j\in \{0, 1, \dotsc, n_m - 8\}$ such that $M_{[j,j+8)} =(A C_A B C_B)^2$.
	Call $J'_m \subseteq J_m$ a subset of maximal cardinality such that
	\begin{equation}\label{eq:estimate_on_Jm}
		\min ((J'_m - J'_m) \cap \NN_{>0}) \ge 8.
	\end{equation}
	We remark that $\# J'_m \ge \frac{1}{8} \# J_m$. Now, if $j\in J'_m$ there
	exists a unique $i = i(j)\in\NN$ such that $n_{i} \in\{j, j+1, j+2\}$ and
	$N_i N_{i+1} \in \{AC_A, BC_B\}$. In particular, if $j, j' \in J'_m$ with
	$j\neq j'$, then $|i(j') - i(j)| \ge 2$, thanks
	to~\eqref{eq:estimate_on_Jm}. Denote $I_m = \{i(j), j\in J'_m\}$, so $\# I_m = \# J'_m$.

	Using \cref{lemma:submultiplicative} recursively
	together with~\eqref{eq:Dnorm_vs_infinitynorm}, \cref{lemma:norm_one} and
	\cref{lemma:AC_ABC_B} we obtain
	\[
		\begin{split}
			\Bigl\|\transpose{M}_{[0,n_m)}\Bigr\|_\infty^{M_{[0, n_m)} \RR^3_{\ge 0}} &=
				\Bigl\|\transpose{N}_{[0,m)}\Bigr\|_\infty^{N_{[0, m)} \RR^3_{\ge 0}} 
				\le 2	\Bigl\|\transpose{N}_{[0,m)}\Bigr\|_D^{N_{[0, m)} \RR^3_{\ge 0}} \\
				&\le 2 \prod_{i\in I_m} \Bigl\|\transpose{(N_i N_{i+1})}\Bigr\|_D^{N_i N_{i+1} \RR^3_{\ge 0}}
				\cdot \prod_{\substack{i\in{0,1,\dotsc, m-1} \\ i\notin I_m, i\notin I_m+1}} \Bigl\|\transpose{N_i}\Bigr\|_D^{N_i \RR^3_{\ge 0}}\\
				&\le 2 \biggl(\frac{4}{5}\biggr)^{\# I_m}.
		\end{split}
	\]

	We can now conclude the proof. Using Birkhoff ergodic theorem, for
	$\mu$-almost every $x\in\{1,2\}^\NN$, we have
	\[
		\lim_{n\to\infty} \frac{1}{n} \sum_{k=0}^{n-8} \chi_{[112211221]} \circ S^k (x) = 
		\lim_{n\to\infty} \frac{1}{n} \sum_{k=0}^{n-1} \chi_{[112211221]} \circ S^k (x) = 
		\mu([112211221]).
	\]
	Hence, for $\mu$-almost every $x\in\{1,2\}^\NN$, and for all $\epsilon>0$,
	there exists an $N$ such that, if $n>N$, then
	\[
		\Biggl| \frac{1}{n} \sum_{k=0}^{n-1} \chi_{[112211221]} \circ S^k (x) - \mu([112211221])\Biggr| < \epsilon,
	\]
	which implies that
	\[
		\begin{split}
			\# J_m &= \sum_{k=0}^{n_m-8} \chi_{[112211221]} \circ S^k (x) \\
			 &\ge \sum_{k=0}^{n-8} \chi_{[112211221]} \circ S^k (x) -1 \\
			 &\ge n(\mu([112211221]) -\epsilon) - 1,
		\end{split}
	\]
	and the proof is complete (we used the definition of the cylinder $[112211221]$ and the coefficients $n_m$ for the first ineuquality).
\end{proof}

The previous work allows us to prove the negativity of the second Lyapunov
exponent for the MCF algorithm. 

\begin{theorem}\label{thm:no_weak_mixing}	
	The MCF algorithm defined by the renormalization algorithm has $\mu$-almost
	everywhere \emph{negative} second Lyapunov exponent. 
\end{theorem}

\begin{proof}
	From the above discussion, using \cref{lemma:measure_cylinder}, for
	$\mu$-almost every word $w$ and every sufficiently small $\epsilon>0$, we
	have that
	\[
	\begin{split}
		\lambda_2 &= \lim_{n\to\infty} \frac{\log{\Bigl\|\transpose{X_n(w)}|_{f^\perp}\Bigr\|}}{n} \\
		&\le \lim_{n\to\infty} \frac{\log{ \Bigl( (n+1) \bigl(\frac{4}{5}\bigr)^{\frac{1}{8}n(\mu([112211221])-\epsilon)-\frac{1}{8}}\Bigr)}}{n} \\
		&= \lim_{n\to\infty} \frac{\log{(n+1)} + \bigl(\frac{1}{8}n(\mu([112211221])-\epsilon)-\frac{1}{8}\bigr)\log{\bigl(\frac{4}{5}\bigr)}}{n} \\
		&< \frac{1}{8}(\mu([112211221])-\epsilon)\log{\biggl(\frac{4}{5}\biggr)}\\
		&< 0,
	\end{split}
	\]
	which proves the statement.

\end{proof}

\section{Hausdorff dimension estimates}\label{sec:HausdorffEstimate}
In this section, we show that the Hausdorff dimension of the gasket defined by
our MCF algorithm is the same as its affinity dimension, and prove estimates for
it. Our approach follows the very recent papers~\cite{Jurga:Gasket,
JiaoLiPanXu:Applications}.

In this section, we will continue to use the simplicial coordinates $(a, b, c)$
to describe the Bruin-Troubetzkoy gasket $\gasket$. We remark again that the
change of coordinates in~\eqref{eq:changeofcoords} on
\cpageref{eq:changeofcoords} given by $(\alpha, \beta) \leftrightarrow (a, b,
c)$ is linear and invertible, hence it is bi-Lipschitz. In particular, the
Hausdorff dimension is preserved.

\subsection{Thermodynamic formalism}\label{sec:thermodynamic_formalism}
In this subsection we recall some general definitions and results
from~\cite{Jurga:Gasket} that we will use. 

Let $\bX$ be a set of matrices in $\SL(3,\RR)$. We say that $\bX$ is
\emph{irreducible} if no proper linear subspace of $\RR^3$ is preserved by all
the matrices in $\bX$.

Given a matrix $X\in\SL(3,\RR)$, let $\alpha_1(X) \ge \alpha_2(X) \ge
\alpha_3(X)$ be its singular values. Then, for $s\ge 0$, the \emph{singular
value function} $\phi^s\colon\SL(3,\RR)\to\RR_+$ is defined as
\begin{equation}\label{eq:singular_value_f}
	\phi^s =
	\begin{cases}
		\Bigl(\frac{\alpha_2(X)}{\alpha_1(X)}\Bigr)^s, & \text{if $0\le s \le 1$},\\	
		\frac{\alpha_2(X)}{\alpha_1(X)}\Bigl(\frac{\alpha_3(X)}{\alpha_1(X)}\Bigr)^{s-1}, & \text{if $1\le s \le 2$},\\
		\Bigl(\frac{\alpha_2(X)\alpha_3(X)}{\alpha_1^2(X)}\Bigr)^{s-1}, & \text{if $s \ge 2$}.
	\end{cases}
\end{equation}

We now define the zeta function $\zeta_\bX\colon [0,\infty)\to[0,\infty]$ by
\[
	\zeta_\bX (s) = \sum_{n=1}^{\infty} \sum_{X\in\bX^n} \phi^s(X).
\]
Finally, the \emph{affinity dimension} $s_\bX$ is the critical exponent of the above
series:
\begin{equation}\label{eq:affinity_dim}
	s_\bX = \inf \{s\ge 0: \zeta_\bX(s)<\infty\}.
\end{equation}

We will use the following result, which uses the strong open set condition
(SOSC), which we will now recall. Given an iterated function system
$\{f_i\}_{i\in\cI}$ on $\RR^2_{>0}$ induced by the projective action of a set of
positive matrices $\bX = \{X_i\}_{i\in\cI} \subset \SL(3,\RR)_{>0}$, satisfies
the \emph{strong open set condition (SOSC)} if there exists an open set
$U\in\RR^2_{>0}$ with non empty intersection with the attractor of the iterated
function systems, denoted $K_\bX$, which satisfies $\cup_{i\in\cI} f_i(U)
\subset U$ and $f_i (U) \cap f_j (U) = \emptyset$ for $i\neq j$ in $\cI$, see,
e.g.,~\cite[Definition~2.5]{Jurga:Gasket} for more details.

\begin{theorem}[\protect{\cite[Theorem~1.3]{Jurga:Gasket}}]\label{thm:Jurga} 
	Suppose a finite set $\bX$ of positive matrices in $\SL(3,\RR)$ generates a 
	semigroup $S_\bX$ which is Zariski dense in $\SL(3,\RR)$ and satisfies the
	SOSC. Then $\dim_H K_\bX = \min \{s_\bX, 2\}$.
\end{theorem}

We say that a finite or countable set of positive matrices $\bX =
\{X_i\}_{i\in\cI}$ is \emph{balanced} if there exists a $c>0$ such that, for all
$i\in\cI$,
\[
	\frac{\min(X_i)_{j,k}}{\max (X_i)_{j,k}} \ge c.
\]
This implies that the singular value function $\phi$ is
\emph{almost-submultiplicative} on the set $\bX$, meaning that there exists a
constant $C<\infty$ such that, for all $A$, $B\in S_\bX$, $\phi(AB) \le C
\phi(A)\phi(B)$, see~\cite[Proposition~2.1]{Jurga:Gasket}.

Another technical condition we will use is the \emph{quasimultiplicativity} of
the singular value function $\phi$ on $S_\bX$. This means that there exists a
finite set of words $W\subset \cI^*$, where $\cI^*$ is the monoid generated by
$\cI$, and a constant $c>0$ such that, for all $p$, $s\in\cI^*$, there exists a
$w\in W$ such that
\[
	\phi^s(X_{pws}) \ge c \phi^s(X_p) \phi^s(X_s).
\]
The irreducibility of $\bX$ implies that the function $\phi^s$ is
quasimultiplicative on $S_\bX$ for all $s\ge 0$,
see~\cite[Section~2.2]{Jurga:Gasket} for details.

Finally, we recall the definition of the pressure $P_\bX\colon[0,\infty)\to\RR$,
given by
	\[
		P_\bX(s) = \lim_{n\to\infty} \frac{1}{n} \log {\Biggl(\sum_{X\in\bX^n} \phi^s(X)\Biggr)},
	\]
where the limit exists since the function $\phi^s$ is almost-submultiplicative.
One can see that $P_\bX$ is a continuous, strictly decreasing, convex function
and that its unique root is exactly the affinity dimension $s_\bX$.

\subsection{The Bruin-Troubetzkoy gasket as the attractor of an IFS}

The graph $\rauzygraph$ in \cref{fig:Rauzy_graph} can be
encoded, using the alphabet $\alphabet = \{A, C_A, B, C_B\}$ by the matrix
\[
	T = (t_{ij}) = \begin{pmatrix}
		1 & 1 & 0 & 0 \\
		0 & 0 & 1 & 1 \\
		0 & 0 & 1 & 1 \\
		1 & 1 & 0 & 0
	\end{pmatrix}.
\]
Considering the paths on the graph $\rauzygraph$ we define a topological Markov
shift, which naturally we call the ARC topological Markov shift. Then, the set
of (one-sided) infinite paths, starting from the vertex $1 = (1, 2, 3)$, on the
graph corresponds to the set
\[
	\sW = \{w = (w_i)_{i=1}^\infty: w_1 = A, C_A, w_i \in\alphabet, t_{w_i w_{i+1}} = 1, \text{ for all } i\ge 1\}.
\]
A word of length $n$, $w = w_1 w_2 \cdots w_n$, is \emph{admissible} if $w_1
= A, C_A$ and 
\[
	t_{w_1w_2} t_{w_2w_3} \cdots t_{w_{n-1}w_n} = 1.
\]
The set of admissible words of length $n$, which corresponds to the set of $n$
length path on the graph $\rauzygraph$, will be denoted $\sW_n$. 

We recall that in \cref{sec:ARC_MCF} we defined the Arnoux-Rauzy-Cassaigne MCF
algorithm. Bearing in mind this algorithm, the Bruin-Troubetzkoy gasket
$\gasket$ is
\[
	\gasket = \bigcup_{w\in\sW} \bigcap_{n=1}^\infty f_{w_1} \circ \cdots \circ f_{w_n} (\Delta).
\]
In other words, it is the \emph{attractor} of the iterated function system
driven by the paths on the graph $\rauzygraph$, using the transformations
in~\eqref{eq:ARC_IFS}.

In the following, we will identity the alphabet $\alphabet$ with the set of
matrices bearing the same names. It is easy to see that the set of matrices
$\alphabet$ is irreducible. Then, we can use the results of the previous
section.

\subsection{Equality between Hausdorff dimension and affinity dimension}
The main result of this section is the following

\begin{theorem}\label{thm:dimH_eq_affdim}
	The Hausdorff dimension of the Bruin-Troubetzkoy gasket $\gasket$ is
	equal to its affinity dimension. That is
	\[
		\dim_H \gasket = s_{\alphabet} = \inf \Biggl\{
		s\ge 0: \sum_{n=1}^{\infty} \sum_{w\in\sW_n} \frac{\alpha_2(w)}{\alpha_1(w)}
			\left(\frac{\alpha_3(w)}{\alpha_1(w)}\right)^{s-1}< \infty
		\Biggr\}.
	\]
\end{theorem}

We will closely follow the strategy used for the analogous result for the Rauzy
gasket in the paper~\cite{Jurga:Gasket}, see her Theorem~1.1.

Let
\[
	\Gamma = \{A^n C_A C_B, C_A B^n C_B, (C_A C_B)^n A\}_{n \ge 1} \subset \SL(3,\RR)
\]
and $\Gamma_N$ the $N$-th truncation of the set:
\[
	\Gamma_N = \{A^n C_A C_B, C_A B^n C_B, (C_A C_B)^n A\}_{1 \le n \le N}.
\]
We will denote by $S_\Gamma$ and $S_{\Gamma_N}$ the semigroups generated by
$\Gamma$ and $\Gamma_N$ respectively. The results recalled in
\cref{sec:thermodynamic_formalism}, applied to these semigroups allow us to
define their affinity dimension as in~\eqref{eq:affinity_dim}. We will denote by
$s_\Gamma$ and $s_{\Gamma_N}$ the affinity dimensions of the semigroups $\Gamma$
and $\Gamma_N$ respectively.

\begin{lemma}\label{lemma:excluding_constants}
	If $K_\Gamma$ denotes the projective limit set of $S_\Gamma$, then
	$\gasket\setminus K_\Gamma$ is countable.
\end{lemma}

\begin{proof}
	The only words on $\alphabet$ that appear in $\sW$ but do not appear as
	combinations of elements in $\Gamma$ are the ones which are eventually
	constantly equal to either $A$, $B$ or $C_A C_B$. Since this set is
	countable, we are done.
\end{proof}

\begin{proposition}\label{prop:simultaneous_conjugacy}
	The matrices in $\Gamma$ can be simultaneously conjugated to a set of
	balanced matrices.
\end{proposition}

\begin{proof}
	We begin by observing that the matrices in $\Gamma$ are non-negative. 
	Consider the matrix
	\[
		M_\epsilon = \begin{pmatrix}
			1 			& -\epsilon & -\epsilon \\
			-\epsilon 	& 1 		& -\epsilon \\
			-\epsilon 	& -\epsilon & 1
		\end{pmatrix},
	\]
	for some sufficiently small $\epsilon$. For instance, $\epsilon \le
	\frac{1}{5}$ is enough.

	A direct computation shows that the entries grow linearly with $n$. This
	implies that we can find two constants $0< c_1 < c_2 < \infty$, which depend
	only on $\Gamma$ and $\epsilon$, such that, for all $X\in\Gamma$ we have
	that $M_\epsilon^{-1} X M_\epsilon = X'$ satisfies $c_1 \le \frac{X'}{n} \le
	c_2$. 
	
	For concreteness, let us compute $A^n C_A C_B M_\epsilon$, we have
	\[
		A^n C_A C_B M_\epsilon = \begin{pmatrix}
			n+1-3n\epsilon  & 2n-(2n+1)\epsilon & n-(3n+1)\epsilon \\
			-\epsilon 		& 1					& -\epsilon \\
			1-2\epsilon		& 1-2\epsilon		& 1-2\epsilon
		\end{pmatrix}.
	\]
	Multiplying by the matrix $M_\epsilon^{-1}$ we see that all the entries,
	once we divided by $n$, are bounded from above and away from $0$. Repeating
	the computation for the matrices $C_A B^n C_B$ and $(C_A C_B)^n A$, we find
	the constants $c_1$ and $c_2$.
\end{proof}

\begin{corollary}\label{cor:sup_affinity_exp}
	We have that $\sup_N s_{\Gamma_N} = s_\Gamma$.
\end{corollary}

\begin{proof}
	Since $\Gamma_N \subseteq \Gamma$, the exponents satisfy $s_{\Gamma_N} \le
	s_\Gamma$. Let $s < s_\Gamma$. By \cref{prop:simultaneous_conjugacy} and the
	discussion after \cref{thm:Jurga}, the function $\phi^s$ is
	almost-submultiplicative on $S_\Gamma$ and $S_N$. Then, $s_\Gamma$ and
	$s_{\Gamma_N}$ are the unique zeros of the respective pressure functions
	$P_\Gamma$ and $P_{\Gamma_N}$. 
	
	Since $\phi^s$ is
	quasimultiplicative,~\cite[Proposition~3.2]{KaenmakiReeve} ensures that
	$0<P_\Gamma (s) = \sup_N P_{\Gamma_N}$, which implies that $P_{\Gamma_N}>0$
	for some $N$.
\end{proof}

\begin{lemma}\label{lemma:equality_affinity_exp}
	We have that $s_\alphabet = s_\Gamma$.
\end{lemma}

\begin{proof}
	Let $S_\alphabet$ be the semigroup generated by $\sW$. Since $\Gamma
	\subseteq S_\alphabet$ we have that $s_\Gamma \le s_\alphabet$. We will now
	show the other direction.
	
	We have
	\[
		\begin{split}
			\zeta_\alphabet (s) &\le \zeta_\Gamma (s) + 
				\sum_{n=1}^{\infty}\sum_{w\in\sW_n} \sum_{k=1}^{\infty} \phi^s(w A^k) + \phi^s(w B^k) + \phi^s(w (C_AC_B)^k) \\
								&\le C\zeta_\Gamma (s) + 
				\sum_{n=1}^{\infty}\sum_{w\in\sW_n} \sum_{k=1}^{\infty} \phi^s(w A^k C_A C_B) + \phi^s(w B^kC_B) + \phi^s(w (C_AC_B)^k A) \\
								&\le 2C\zeta_\Gamma(s),
		\end{split}	
	\]
	where we used that we can find a constant $C < \infty$, that only depend on
	the matrices in $\alphabet$ and $s$, such that for all matrices
	$X\in\alphabet$ and $Y\in\SL(3,\RR)$, we have $\phi^s(Y) \le C\phi^s(YX)$,
	see, e.g.,~\cite[Lemma~1]{BochiGourmelon}. The above inequalities imply
	$s_\alphabet \le s_\Gamma$ and we are done.
\end{proof}

\begin{proposition}\label{prop:Zariski_density}
	For all sufficiently large $N$, the subgroup $S_N$ generated by $\Gamma_N$
	is Zariski dense in $\SL(3,\RR)$.
\end{proposition}

\begin{proof}
	We begin by recalling that the Zariski closure of any subgroup is an
	algebraic group. Let $G$ be the Zariski closure of $S_\Gamma$, and
	$\mathfrak{g}$ its Lie algebra, which corresponds to the tangent space to
	the identity. It is clear that $\mathfrak{g} \subseteq \sl(3, \RR)$, where
	$\sl(3, \RR)$ is the Lie algebra of the Lie group $\SL(3,\RR)$ is the
	$8$-dimensional algebra of the $3\times 3$ matrices with zero trace and the
	usual matrix conmutador as Lie bracket: $[X,Y] = XY - YX$. We will show that
	$\mathfrak{g} = \sl(3, \RR)$, by finding $8$ linearly independent matrices
	in $\mathfrak{g}$.

	Let us consider the matrices
	\[
		A^n (C_A C_B) = \begin{pmatrix}
			n+1 & 2n & n \\
			0 & 1 & 1 \\
			1 & 1 & 1
		\end{pmatrix},
	\]
	for any $n\in\NN$. Let $P$ be a polynomial that is zero on all the points of
	$S_\Gamma$. Then we can form the real polynomial
	\[
		q(x) = P \left(\begin{pmatrix}
			x+1 & 2x & x \\
			0 & 1 & 1 \\
			1 & 1 & 1
		\end{pmatrix}\right).
	\]
	Since $q(n)=0$ for all $n\in\NN$, $q\equiv 0$, which implies that the
	matrices
	\[
		\gamma(x) = 
		\begin{pmatrix}
			x+1 & 2x & x \\
			0 & 1 & 1 \\
			1 & 1 & 1
		\end{pmatrix}
	\]
	form a curve inside the real algebraic group $G$. Then
	\[
		X_1 = \frac{d}{dx} \gamma(0)^{-1}\gamma(x)\biggr|_{x=0} = \begin{pmatrix}
			1  & 2  & 1 \\
			0  & 0  & 0 \\
			-1 & -2 & -1
		\end{pmatrix}\in\mathfrak{g}.
	\]
	Similarly, by considering $C_A B^n C_B$ and $(C_A C_B)^n$ respectively, we
	obtain that
	\[
		X_2 = \begin{pmatrix}
			0  & 0  & 0 \\
			1  & 1  & 1 \\
			-1 & -1 & -1
		\end{pmatrix}
		\qquad
		\text{and}
		\qquad
		X_3 = \begin{pmatrix}
			0 & 0 & 0 \\
			0 & 0 & 0 \\
			1 & 1 & 0
		\end{pmatrix}
	\]
	are in $\mathfrak{g}$.
	We now consider the commutators
	\[
		\begin{split}
			X_4 &= [X_1, X_2] = \begin{pmatrix}
				1  & 1  & 1 \\
				0  & 0  & 0 \\
				-1 & -1 & -1
			\end{pmatrix},
			\qquad
			X_5 = [X_1, X_3] = \begin{pmatrix}
				1  & 1  & 0 \\
				0  & 0  & 0 \\
				-2 & -3 & -1
			\end{pmatrix},\\
			X_6 &= [X_2, X_3] = \begin{pmatrix}
				0  & 0  & 0 \\
				1  & 1  & 0 \\
				-2 & -2 & -1
			\end{pmatrix},
			\qquad
			X_7 = [X_3, X_4] = \begin{pmatrix}
				-1 & -1 & 0 \\
				0  & 0  & 0 \\
				-2 & -2 & 1
			\end{pmatrix},\\
			X_8 &= [X_2, X_5] = \begin{pmatrix}
				-1 & -1 & -1 \\
				-1 & -2 & 0 \\
				4  & 5  & 3
			\end{pmatrix}.
		\end{split}
	\]

	It can be checked that the set $\{X_i\}_{i=1}^8$ is a linearly independent
	subset of $\mathfrak{g}$, and hence $\mathfrak{g}=\sl(3,\RR)$. Thus,
	$S_\Gamma$ is Zariski dense inside $\SL(3,\RR)$.

	To conclude the proof, we remark that, since $S_\Gamma$ is a subsemigroup of
	$S_\alphabet$, the latter is also Zariski dense inside $\SL(3,\RR)$. Since
	$\SL(3,\RR)$ is a (Zariski) closed and connected subgroup of $\GL(3,\RR)$,
	density of $\Gamma_N$ for sufficiently large $N$ follows
	from~\cite[Lemma~3.7]{MorrisSert:Variational}.
\end{proof}

We can now prove the main result of this section.

\begin{proof}[Proof of \cref{thm:dimH_eq_affdim}]
	From \cref{prop:simultaneous_conjugacy} and \cref{prop:Zariski_density}, for
	sufficiently large $N$ we have can simultaneously conjugate every $\Gamma_N$
	to a subset of \emph{positive} matrices in $\SL(3,\RR)$ which satisfies the
	SOSC and that generate a Zariski dense subgroup of $\SL(3,\RR)$. Then,
	\cref{thm:Jurga}, together with \cref{cor:sup_affinity_exp} and
	\cref{lemma:equality_affinity_exp} yield
	\[
		\dim_H \gasket \ge \sup_N \dim_H K_{\Gamma_N} = \sup_N s_{\Gamma_N} = s_\Gamma = s_\alphabet.
	\]

	Let us show the reverse inequality. Since $\gasket \setminus K_\Gamma$ is
	countable, we have that $\dim_H \gasket = \dim_H K_\Gamma$. By
	\cref{prop:simultaneous_conjugacy} we can simultaneously conjugate $\Gamma$
	to $\Gamma_\epsilon$, a set of positive matrices in $\SL(3,\RR)$. These
	matrices send the positive cone into a compact subset of itself, so one can
	reason as in~\cite[Section~6.1.2]{Jurga:Gasket} to show that $\dim_H
	K_\Gamma \le s_\Gamma$. Finally, by \cref{lemma:equality_affinity_exp}, we
	have $s_\Gamma = s_\alphabet$. Then, $\dim_H \gasket \le s_\alphabet$.
\end{proof}

\subsection{A lower bound on $\dim_H \gasket$.}
We recall that we denote by $\cE = \{e_1, e_2, e_3\}$ the standard base of
$\RR^3$. The corresponding elements in $\PP(\RR^3)$ will be denoted by $E_i =
\RR e_i$. In this section, it will be more convenient to use a different set of
generators for the semigroup $S_\Gamma$. Let 
\[
	D_1 = A, \quad D_2^n = C_A B^n C_B, \quad D_3 = C_A C_B = \begin{pmatrix}
		1 & 0 & 0 \\
		0 & 1 & 0 \\
		1 & 1 & 1 
	\end{pmatrix},
\]
for any $n \ge 1$. Then, $\{D_1, D_2^n, D_3, n\in\NN\}$ still generates the
semigroup $s_\Gamma$. Since we will also use their transpose, we list them,
to help the reader:
\[
	\transpose{D_1} = \begin{pmatrix}
		1 & 0 & 0 \\
		1 & 1 & 0 \\
		1 & 0 & 1
	\end{pmatrix},
	\quad
	\transpose{D_2^n} = \begin{pmatrix}
		1 & n & 1 \\
		0 & n+1 & 1 \\
		0 & n & 1
	\end{pmatrix},
	\quad
	\transpose{D_3} = \begin{pmatrix}
		1 & 0 & 1 \\
		0 & 1 & 1 \\
		0 & 0 & 1
	\end{pmatrix}.
\]

The following result follows from direct computations and will be left to the
reader.

\begin{lemma}\label{lemma:deltaprime}
	The matrices $\{D_1, D_2^n, D_3\}$ and their transpose preserve the simplex
	$\Delta$. Moreover, $\Bigl\{\transpose{D_1}, \transpose{D_2^n},
	\transpose{D_3}\Bigr\}$ also preserves the open sub-simplex $\Delta'$ with
	vertices $(0:1:1)$, $(1:0:1)$ and $(1:1:0)$.
\end{lemma}

From this, we obtain the following useful corollary. Let us introduce some
notation we will need. Let $\gamma = \gamma_1 \gamma_2 \cdots \gamma_k \in
S_\Gamma$, then $\transpose{\gamma} = \transpose{\gamma_k} \cdots
\transpose{\gamma_2} \transpose{\gamma_1}$. Moreover, given $1\le m \le k$, we
denote $\transpose{\gamma}_{[1,m)} = \transpose{\gamma_k} \cdots
\transpose{\gamma_{k-m+1}}$. We stress that, in the previous notation, we first
take the transpose and then cut the product after $m$ terms. 

\begin{corollary}\label{cor:insidedeltaprime}
	For any $i=1,2,3$, for any $k \ge 1$ and any $\gamma\in S_\Gamma$ of length 
	$k$, if $\gamma_j = i$ for some $j$, then $\transpose{\gamma} E_i \in \Delta'$.
\end{corollary}

\begin{proof}
	A direct computation shows that
	\[
		\transpose{D_1} E_1 = \transpose{D_3} E_3 = (1:1:1)\in\Delta',
	\]
	and
	\[
		\transpose{D_2^n} E_2 = (n:n+1:n)\in\Delta',
	\]
	as we wanted.
\end{proof}

\begin{proposition}\label{prop:gammainDelta}
	Let $\gamma\in S_\Gamma$, and assume that its last $m$ letters are
	not the same. Then, $\transpose{\gamma}\gamma \Delta \subset \transpose{\gamma}_{[1,m)}
	\Delta \cap \overline{\Delta'}$.
\end{proposition}

\begin{proof}
	The inclusion $\transpose{\gamma}\gamma \Delta \subset \transpose{\gamma}_{[1,m)} \Delta$
	holds trivially. Hence, we need to show that $\transpose{\gamma}\gamma \Delta \subset
	\Delta'$. If all the letters appear in $\gamma$, then we can conclude by
	\cref{cor:insidedeltaprime}.
	
	Since not all the last $m$ letters of $\gamma$ are the same, we must have
	that two among $\{D_1, D_2^n, D_3\}$ appear in the last $m$ digits of
	$\gamma$. We consider each case separately,

	\textbf{Case 1: only $D_1$ and $D_2^n$ appear.}
	In this case, we actually have to distinguish whether $D_1$ or $D_2^n$
	occurs first. If we have $D_1 D_2^n$ inside $\gamma$, then
	$\transpose{\gamma}$ contains $\transpose{D_2^n} \transpose{D_1}$. Since
	\[
		\transpose{D_2^n} \transpose{D_1} = \begin{pmatrix}
			1 & n & 1 \\
			0 & n+1 & 1 \\
			0 & n & 1
		\end{pmatrix}
		\begin{pmatrix}
			1 & 0 & 0 \\
			1 & 1 & 0 \\
			1 & 0 & 1
		\end{pmatrix}
		=
		\begin{pmatrix}
			n+2 & n & 1 \\
			n+2 & n+1 & 1 \\
			n+1 & n & 1
		\end{pmatrix},
	\]
	we have that $\transpose{D_2^n} \transpose{D_1} E_i \in \Delta'$ for
	$i=1,2,3$, as we wanted.

	Similarly, if $D_2 D_1^n$ is contained inside $\gamma$, then
	$\transpose{\gamma}$ contains $\transpose{D_1} \transpose{D_2^n}$. Since
	\[
		\transpose{D_1} \transpose{D_2^n} = \begin{pmatrix}
			1 & n & 1 \\
			1 & 2n+1 & 2 \\
			1 & 2n & 2
		\end{pmatrix},
	\]
	we have that $\transpose{D_1} \transpose{D_2^n} E_i \in \Delta'$ for
	$i=1,2,3$, and we are done with this case.

	\textbf{Case 2: only $D_3$ and $D_2^n$ appear.}
	Let us define the set
	\[
		\nabla_x = \{(x:y:z) \in \Delta: x \le y+z\}.
	\]
	It can be checked that this set is invariant under the action of $D_2^n$,
	$D_3$ and their transposes. Moreover, $D_2^n \Delta \subset \nabla_x$ and
	$D_3 \Delta \subset \nabla_x$. Hence, $\transpose{\gamma}\gamma \Delta
	\subset \nabla_x$. Now we compute
	\[
		\transpose{D_3} \transpose{D_2^n} = \begin{pmatrix}
			1 & 2n & 2 \\
			0 & 2n+1 & 2 \\
			0 & n & 1
		\end{pmatrix}
		\qquad
		\text{and}
		\qquad
		\transpose{D_2^n} \transpose{D_3} = \begin{pmatrix}
			1 & n & n+2 \\
			0 & n+1 & n+2 \\
			0 & n & n+1
		\end{pmatrix}.
	\]
	We observe that, in both cases, the last two columns belong to $\Delta'$,
	while the first coordinate is invariant. In other words,
	$\transpose{\gamma}_{[1,m)} E_2$, $\transpose{\gamma}_{[1,m)} E_3 \in
	\Delta'$, whereas $\transpose{\gamma}_{[1,m)} E_1 = E_1$. Finally we have
	$\transpose{\gamma}\gamma \Delta \subset \transpose{\gamma}_{[1,m)} \Delta
	\cap \nabla_x = \transpose{\gamma}_{[1,m)} \Delta \cap \overline{\Delta'}$,
	and we are done.

	\textbf{Case 3: only $D_1$ and $D_2$ appear.}
	This case can be treated as the previous one, replacing $D_2^n$ by $D_1$ and
	$\nabla_x$ by $\nabla_z$, which is defined analogously.
\end{proof}

The following is the key technical result of this section.

\begin{lemma}\label{lemma:lastdigitsgamma}
	For every $m\in\NN$, there exists an $\epsilon_m>0$ such that, for all
	$\gamma\in\Gamma$, if the last $m$ letters of $\gamma$ are not the same,
	then we have
	\[
		\| \gamma e_i \| \ge \epsilon_m \alpha_1(\gamma),
	\] 
	for $i=1,2,3$.
\end{lemma}

\begin{proof}
	Using the $KA^+K$ Cartan decomposition of $\SL(3,\RR)$, we can write every
	matrix $X\in\SL(3,\RR)$ as $\tilde{k}_X a_X k_X$ where $\tilde{k}_X$,
	$k_X\in\SO(3,\RR)$ and $a_X$ is the diagonal matrix made by the singular
	values $\alpha_1(X)\ge\alpha_2(X)\ge\alpha_3(X)$.
	By~\cite[Lemma~14.2]{BenoistQuint:RandomBook}, one has
	\[
		\|\gamma e_i \| \ge \|\gamma\| d(E_i, H_\gamma),
	\]
	where $H_\gamma = k_\gamma^{-1} (E_i^\perp)$ is a repelling
	hyperplane for $\gamma$ and
	\[
		d(\RR v, \RR w) = \frac{\|v \wedge w \|}{\|v\| \|w\|},
	\]
	with $\|\cdot\|$ the standard Euclidean norm on $\RR^3$ and the induced one
	on $\wedge^2 \RR^3$.

	One can check that $(\tilde{k}_\gamma E_i)^\perp =
	(V_{\transpose{\gamma}})^\perp = H_\gamma$. Hence, it is enough to
	check that the angle between $E_i$ and $V_{\transpose{\gamma}}$ is
	bounded away from $\frac{\pi}{2}$. Since $E_i^\perp$ is the span of $E_j$,
	for $j\neq i$, which is an edge of the simplex $\Delta$, it is enough to
	show that $d(V_{\transpose{\gamma}}, \partial\Delta)$ is bounded from
	below by a constant that only depends on $m$, not on $\gamma$. By
	definition, $V_{\transpose{\gamma}}$, is the attracting fixed point of
	$\transpose{\gamma}\gamma$. In particular,
	$V_{\transpose{\gamma}}\in\transpose{\gamma}\gamma\Delta$.

	By \cref{prop:gammainDelta}, we have that $\transpose{\gamma}\gamma \Delta
	\subset \transpose{\gamma}_{[1,m)} \Delta \cap \overline{\Delta'}$. We
	remark that $\transpose{\gamma}_{[1,m)} \Delta \cap \overline{\Delta'}$ is a
	quadrilateral which does not intersect the boundary of the simplex. Thus,
	$d(\partial\Delta, \transpose{\gamma}_{[1,m)} \Delta \cap
	\overline{\Delta'}) > 0$. Since, for any given $m$, there exists only
	finitely many $\gamma$ of length $m$, we can find a $d_m>0$ such that
	\[
		d(V_{\transpose{\gamma}}, \partial\Delta) \ge d\Bigl(\transpose{\gamma}\gamma\Delta, \partial\Delta\Bigr)
			\ge d\Bigl(\transpose{\gamma}_{[1,m)} \Delta \cap \overline{\Delta'}, \partial\Delta\Bigr) > d_m,
	\]
	which completes the proof.
\end{proof}

We will need an estimation of the distortion of the simplex $\Delta$ by an
element $\gamma\in\Gamma$.

\begin{lemma}\label{lemma:Distortion}
	Assume that the last two letters of $\gamma$ are not the same. Then, there
	exists a constant $C_2 > 1$ such that:
	\begin{enumerate}
		\item $\diam(\gamma\Delta) \le C_2
		\frac{\alpha_2(\gamma)}{\alpha_1(\gamma)}$.
		\item $\area(\gamma\Delta) \le C_2 \alpha_1(\gamma)^{-3}$.
	\end{enumerate}
\end{lemma}

\begin{proof}
	We begin with the first point. It is enough to show that $d(\gamma E_i,
	\gamma E_j) \le C_2 \frac{\alpha_2(\gamma)}{\alpha_1(\gamma)}$, for
	$i, j = 1, 2, 3$.
	We have that
	\[
		\begin{split}
			d(\gamma E_i, \gamma E_j) &= \frac{\|\gamma e_i \wedge \gamma e_j\|}{\|\gamma e_i\| \|\gamma e_j\|} \\
				&\le \frac{\alpha_1(\gamma)\alpha_2(\gamma) \|e_i \wedge e_j \|}{\|\gamma e_i\| \|\gamma e_j\|} \\
				&\le \frac{\alpha_1(\gamma)\alpha_2(\gamma)}{\epsilon_2^2 (\alpha_1(\gamma))^2} \\
				&= \epsilon_2^{-2}\frac{\alpha_2(\gamma)}{\alpha_1(\gamma)},
		\end{split}
	\]
	where we used \cref{lemma:lastdigitsgamma} in the second inequality.

	To prove the second point, we begin with the following elementary
	geometrical fact. Let $x,y,z\in\RR^3\setminus\{0\}$, then the area of the
	triangle $\stackrel{\triangle}{xyz}$ with vertices $x$, $y$ and $z$ is given
	by
	\[
		\area(\stackrel{\triangle}{xyz}) = \frac{\|x \wedge y \wedge z\|}{2d_E(0,\stackrel{\triangle}{xyz})},
	\]
	where $d_E$ is the distance from the origin to the plane containing the
	three points. Slightly abusing the notation, we identify the projective
	simplex $\Delta$ with the ordinary $3$-simplex in $\RR^3$. Let $x_\gamma =
	\gamma e_1$, $y_\gamma = \gamma e_2$ and $z_\gamma = \gamma e_3$ the three
	vertices of $\gamma\Delta$. By definition,
	\[
		x_\gamma = \frac{\gamma e_1}{\|\gamma e_1\|_1},
	\]
	where $\|\cdot\|_1$ is the $\ell^1$-norm on $\RR^3$, and similarly for the
	other points. Hence
	\[
		\|x_\gamma \wedge y_\gamma \wedge z_\gamma \| = \frac{\|\gamma e_1 \wedge \gamma e_2 \wedge \gamma e_3\|}{\prod_{i=1}^3\|\gamma e_i \|_1}
	\]
	Since $\|\gamma e_1 \wedge \gamma e_2 \wedge \gamma e_3\| = \|e_1 \wedge e_2
	\wedge e_3\| = 1$, and the entries of $\gamma$ are non-negative, then
	$\|\gamma e_i \|_1 \ge \|\gamma e_i\|$, so by \cref{lemma:lastdigitsgamma}
	we are done.
\end{proof}

We will need the following geometrical result, which follows from the definition
of Hausdorff dimension and was proven in~\cite[Lemma~4.1]{PollicottSewell}.

\begin{lemma}\label{lemma:PS} For every $\delta > 0$, there exists a $C_\delta >
	0$ such that, for all $\gamma\in\Gamma$, there exists a finite open cover
	$\{B_i(\gamma)\}_{i=1,\dotsc,k}$ of $\gamma\Delta$ with
	$\diam(B_i(\gamma))\le\diam(\gamma\Delta)$ such that
	\[
		\sum_{i=1}^k \diam^{1+\delta} B_i(\gamma) \le 
			c_\delta \cdot \diam^{1-\delta} \gamma\Delta \cdot \area^\delta \gamma\Delta.
	\]
\end{lemma}

Exactly as in \cref{lemma:excluding_constants}, we can decompose $\gasket$ into
a set of \emph{nice} points $\gasket_{\text{nice}}$ whose coding is not
eventually constant and a countable set. However, we remark that
$\gasket_{\text{nice}} \neq K_\Gamma$, since we have switched the generating
set.

Let
\[
	\Gamma^m = \{\gamma\in\Gamma: \text{ the last two digits of $\gamma$ are different and } \diam\gamma\Delta \le 1/m\},
\]
and consider the two families of coverings:
\[
	\mathcal{U}_m = \{B_i(\gamma), \gamma\in\Gamma^m\}, 
	\quad \text{and} \quad
	\mathcal{U}'_m = \{\gamma\Delta, \gamma\in\Gamma^m\},
\]
with $B_i(\gamma)$ defined by \cref{lemma:PS}. We define
\[
	Y = \bigcap_{m=1}^\infty \bigcup_{U\in\mathcal{U}_m} U,
	\quad \text{and} \quad
	Y' = \bigcap_{m=1}^\infty \bigcup_{U\in\mathcal{U}'_m} U.
\]	
Then $\mathcal{U}_m$ is a Vitali cover of $Y$: for every $y\in Y$ and every
$\delta>0$, there exists some $U\in\mathcal{U}_m$ such that $\diam U < \delta$
and $y\in U$. Similarly, $\mathcal{U}'_m$ is a Vitali cover of $Y'$. Moreover,
by construction $Y\supset Y'$.

We have

\begin{lemma}\label{lemma:nice_points}
	We have the inclusion $\gasket_{\text{nice}} \subset Y$.
\end{lemma}

\begin{proof}
	Let $x$ be a point in $\gasket_{\text{nice}}$. Then, its coding with respect
	to $\{D_1, D_2, D_3\}$ is not eventually constant. In particular, there are
	infinitely many pairs of adjacent letters in its coding which are different
	one from the other. By dividing the coding into subwords after each of these
	pairs appears, we form infinitely many words $\gamma\in\Gamma$ whose two
	last letters are different. Thus $x\in\gamma\Delta$. Moreover, since
	$\diam\gamma\Delta\to 0$ as the length of $\gamma$ increases, we see that
	$x\in Y' \subset Y$.
\end{proof}

Let $\Gamma_0 = \langle D_1, D_3 \rangle$ be the semigroup generated by the
matrices $D_1 = A$ and $D_3 = C_A C_B$. Then, the arc $I = I(E_1, E_3) = \{\RR(s
e_1 + t e_3), s,t\in\RR_{\ge 0}\}$ is preserved by $\Gamma_0$.

\begin{lemma}
	There exists an $\epsilon > 0$ such that, for all $\gamma\in\Gamma_0$ having
	the last two digits different from each other, we have
	\[
		\alpha_2(\gamma) \ge \epsilon, 
		\qquad \text{and} \qquad
		\epsilon |\gamma I| \le \frac{1}{\alpha_1(\gamma)^2},
	\]
	where $|\gamma I|$ is the length of the arc $\gamma I$.
\end{lemma}

\begin{proof}
	Since the matrices $D_1$ and $D_3$ preserve $I$ and their restriction to the
	subspace generated by $E_1$ and $E_3$ has determinant one, by
	\cref{lemma:lastdigitsgamma}, we have
	\[
		|\gamma I| = d(\gamma E_1, \gamma E_3) = 
			\frac{\|\gamma e_1 \wedge \gamma e_3\|}{\|\gamma e_1\| \|\gamma e_3\|} 
			= \frac{1}{\|\gamma e_1\| \|\gamma e_3\|} \le \frac{1}{\alpha_1(\gamma)^2}.
	\]
	We recall that the top singular value gives the operator norm of a matrix.
	Hence, $\alpha_1 (\gamma) \ge \|\gamma e_i \|_1$ for $i=1$, $2$, $3$. So, as
	in the proof of \cref{lemma:Distortion}, we obtain that
	\[
		\area(\gamma\Delta) = C_2 \frac{\|\gamma e_1 \wedge \gamma e_2 \wedge \gamma e_3\|}{\prod_{i=1}^{3}\|\gamma e_i\|_1} 
			\ge C_2 \frac{1}{\alpha_1(\gamma)^3}.
	\]
	Combining the last two inequalities, we obtain
	\[
		\max\{|\gamma I(E_1, E_2)|, |\gamma I(E_2, E_3)|\} \ge \frac{\area(\gamma\Delta)}{|\gamma I|} \ge C_2 \frac{1}{\alpha_1(\gamma)}.
	\]
	Finally, using the first part of \cref{lemma:Distortion}, we have
	\[
		\alpha_2(\gamma) \ge C_2 \alpha_1 (\gamma )\diam(\gamma\Delta) 
			\ge C_2 \alpha_1 (\gamma ) \max\{|\gamma I(E_1, E_2)|, |\gamma I(E_2, E_3)|\}
			\ge C_2^2 \frac{1}{\alpha_1(\gamma)}
	\]
	and we are done.
\end{proof}

We are now ready to conclude and prove the main result of this section.

\begin{theorem}
	The Hausdorff dimension of the Bruin-Troubetzkoy gasket $\gasket$ is
	greater than $3/2$:
	\[
		\dim_H \gasket = s_\alphabet \ge \frac{3}{2}.
	\]
\end{theorem}

\begin{proof}
	Let $\gamma\in\Gamma_0$ be an element whose last two digits are different,
	the singular value function~\eqref{eq:singular_value_f} we have
	\[
		\phi^{3/2} (\gamma) = \frac{\alpha_2(\gamma)}{\alpha_1(\gamma)} \biggl(\frac{\alpha_3(\gamma)}{\alpha_1(\gamma)}\biggr)^{1/2}
		= \frac{\alpha_2(\gamma)^{1/2}}{\alpha_1(\gamma)^2} \ge \frac{\epsilon^{1/2}}{\alpha_1(\gamma)^2} \ge \epsilon^2 |\gamma I|.
	\]
	Then
	\begin{multline*}
		\sum_{n=1}^{\infty} \sum_{w\in\sW_n} \phi^{3/2}(w) 
		\ge \sum_{\gamma\in\Gamma_0} \phi^{3/2}
		\ge \sum_{\substack{\gamma\in\Gamma_0 \\ \text{last two digits different}}} \phi^{3/2} \\
		\ge \epsilon^2 \sum_{\substack{\gamma\in\Gamma_0 \\ \text{last two digits different}}} |\gamma I|.
	\end{multline*}

	We remark that, by \cref{lemma:nice_points}, every point inside $\gamma \cap
	\gasket_{\text{nice}}$ is covered infinitely many times by
	\[
		\{\gamma I, \gamma\in\Gamma_0 \text{ with the last two digits different from each other}\}.
	\]
	Since $I \setminus (I\cap \gasket_{\text{nice}})$ is countable, the series
	$\ge \sum_{\gamma\in\Gamma_0} \phi^{3/2}$ diverges and $\dim_H \gasket \ge
	\frac{3}{2}$, as we wanted.
\end{proof}

\printbibliography
\end{document}

\typeout{get arXiv to do 4 passes: Label(s) may have changed. Rerun}